\documentclass{article}
\usepackage{amsmath,amssymb,enumerate}
\usepackage{graphicx}
\usepackage{float}
\usepackage[latin1]{inputenc}
\usepackage[T1]{fontenc}
\usepackage[compatibility=false]{caption}
\usepackage{subcaption}
\usepackage{dcolumn}
\usepackage{color, colortbl}
\usepackage{geometry}
\usepackage{lmodern}
\usepackage{hyperref}
\usepackage{slashbox}
\usepackage[para,online,flushleft]{threeparttable}
\usepackage{multicol}

\newcommand{\fX}{\ensuremath{f_X}}

\newcommand{\fbeta}{\ensuremath{f_{\beta}}}
\newcommand{\ind}{\ensuremath{\mathbf{1}}}
\newcommand{\cH}{\ensuremath{\mathcal{H}}}

\newcommand{\hatbeta}{\ensuremath{\widehat{\beta}}}
\newcommand{\hatfbeta}{\ensuremath{\widehat{\fbeta}}}
\newcommand{\Esp}{\ensuremath{\mathbb{E}}}
\newcommand{\Prob}{\ensuremath{\mathbb{P}}}

\newcommand{\bR}{\ensuremath{\mathbb{R}}}
\newcommand{\Quad}{\ensuremath{\Xi}}

\newcommand{\ConstQuadNb}{\ensuremath{C_\Quad}}
\newcommand{\ConstNeedNorm}{\ensuremath{C}}

\newcommand{\ConstNeedEquiv}{\ensuremath{C'}}
\newcommand{\ConstNeedEquivTwo}{\ensuremath{C''}}

\newcommand{\infzr}{\ensuremath{(z\wedge r)}}

\newcommand{\harm}{\ensuremath{h}}

\newcommand{\fbetam}{\ensuremath{{\fbeta^{-}}}}

\newcommand{\hatfbetam}{\ensuremath{\widehat{\fbetam}}}

\newcommand{\odd}{\text{ odd}}
\newcommand{\even}{\text{ even}}

\newcommand{\E}{\mathbb{E}}

\newcommand{\ConstTsigmap}{C_{\sigma,p}}
\newcommand\myeq{\mathrel{\stackrel{\makebox[0pt]{\mbox{\normalfont\tiny def}}}{=}}}
\def\xL{{\rm L}}

\def\xR{{\mathbb R}}

\def\xN{{\mathbb N}}

\def\supp{{\rm supp}}
\def\xS1{{\mathbb S}^1}
\def\xSd{{\mathbb S}^{d-1}}
\def\xRd{{\mathbb R}^{d}}

\def\xLtwo{{{\rm L}^2}}

\def\xLn{{\rm L}}
\def\xCn{{\rm C}}
\def\xWn{{\rm W}}

\def\xLinfty{{\rm L}^{\infty}}
\newtheorem{theorem}{Theorem}
\newtheorem{proposition}[theorem]{Proposition}
\newtheorem{lemma}[theorem]{Lemma}

\newtheorem{definition}[theorem]{Definition}

\newtheorem{assumption}[theorem]{Assumption}

\title{Adaptive estimation in the nonparametric random coefficients
binary choice model by needlet thresholding}

\author{Eric
  Gautier\thanks{Toulouse School of Economics, Universit\'e Toulouse
    Capitole, 21 all\'ee de Brienne, 31015 Toulouse, France,
    \url{Eric.Gautier@tse-fr.eu}}\ \thanks{Eric Gautier acknowledges
    financial support from the grants ERC POEMH and ANR-13-BSH1-0004}\ 
\thanks{The authors are grateful to Andrii Babii (TSE) for research assistance on the simulations.}
\ and
Erwan Le Pennec\thanks{CMAP / D\'epartement de Math\'ematiques Appliqu\'ees, Ecole Polytechnique,
91128 Palaiseau, France.
\url{Erwan.Le-Pennec@polytechnique.edu}}
\ \footnotemark[3]
}

\date{November 2017}

\begin{document}

\maketitle

\begin{abstract}
In the random coefficients binary choice model, a binary variable equals 1 iff an index $X^\top\beta$  is positive.
The vectors $X$ and $\beta$ are independent and belong to the sphere $\mathbb{S}^{d-1}$ in $\mathbb{R}^{d}$.
We prove lower bounds on the minimax risk for estimation of the density $f_{\beta}$ over Besov bodies where the loss is a power of the $\xL^p(\mathbb{S}^{d-1})$ norm for $1\le p\le \infty$.  
We show that a hard thresholding estimator based on a needlet expansion with data-driven thresholds achieves these lower bounds up to logarithmic factors.
\end{abstract}

\tableofcontents

\section{Introduction}\label{s1}
Discrete choice models (see, \emph{e.g.}, \cite{MF}) have applications in many areas
ranging from planning of public transportation, economics of industrial organizations, evaluation of
public policies, among others.
This paper considers the binary choice model. There, agents (consumer, firm, country, etc.) choose between two exclusive alternatives 1 or -1 (\emph{e.g.}, buying a good or not) the one that yields the highest utility. 
The utility that an agent $i$ gets from choosing alternative -1 (resp. from choosing 1) is assumed to have the form
\begin{align}\label{eq:utility}
u_{-1,i}=z_{-1,i}^\top\gamma_i+\epsilon_{-1,i} \quad\mathrm{(resp.\  }
u_{1,i}=z_{1,i}^\top\gamma_i+\epsilon_{1,i}),
\end{align}
where $z_{-1,i}$ (resp. $z_{1,i}$) is a vector of $d-1$ characteristics of alternative -1 (resp. 1) for agent $i$, $d\ge2$, $\gamma_i$ are preferences of agent $i$ for the characteristics, and $\epsilon_{-1,i}$ and $\epsilon_{1,i}$ absorb both the usual error terms and constants. In \eqref{eq:utility}, the preferences are allowed to vary across individuals; namely, they are heterogeneous. This translates into a vector of coefficients $\gamma$ indexed by $i$ that we assume random.
The characteristics of the alternatives are indexed by the agents, for example they can be characteristics of two goods that a consumer has to choose upon interacted with individual characteristics like age or distance. 
We assume that the random coefficients and errors are independent from the characteristics.
The statistician observes a sample of characteristics and choices for agents $i=1,\dots,n$, but $\gamma_i$, $u_{1,i}$, and $u_{-1,i}$ are not observed.  Observing the choices corresponds to observing the sign $y_i$ of the net utility $u_{1,i}-u_{-1,i}$. Indeed, agent $i$ prefers 1 ($y_i=1$) if and only if the net utility for 1 is positive, \emph{i.e.},
\begin{equation}\label{eqnet}
u_{1,i}-u_{-1,i}=\epsilon_{1,i}-\epsilon_{-1,i}+(z_{1,i}-z_{-1,i})^\top\gamma_i>0,
\end{equation} 
and prefers -1 ($y_i=-1$) when $$u_{1,i}-u_{-1,i}<0.$$ 
We assume that the probability that $\left|(\epsilon_{1,i}-\epsilon_{-1,i}, \gamma_i^\top)^\top\right|$ is the 0 and thus that  agent $i$ is indifferent (\emph{i.e.}, $u_{1,i}-u_{-1,i}=0$) on a set of 0 probability. 
Hence, the linear random coefficients binary choice model is
\begin{equation}\label{model}
y_i= {\rm sign} \left(x_i^\top\beta_i\right),
\end{equation}
where, for a real number $a$, ${\rm sign}(a)$ is 1 if $a>0$, -1 if $a<0$, and is 0 if $a=0$,
\begin{align*}
x_i&=(1, (z_{1,i}-z_{-1,i})^\top)^\top/\left|(1, (z_{1,i}-z_{-1,i})^\top)^\top\right|,\\ 
\beta_i&=(\epsilon_{1,i}-\epsilon_{-1,i}, \gamma_i^\top)^\top/\left|(\epsilon_{1,i}-\epsilon_{-1,i}, \gamma_i^\top)^\top\right|,
\end{align*}  
and $|\cdot|$ is the Euclidean norm in $\xRd$. 
Like in \cite{BFH,BM,GK,IT} among others,
we consider a nonparametric specification of the joint distribution of $\beta$ and this model is more general than the Logit, Probit, and Mixed-Logit models. Note that it is important to avoid restricting the dependence between the coordinates of $(\epsilon_{1}-\epsilon_{-1},\gamma^\top)$ since they can be functions of a deep heterogeneity parameter (\emph{e.g.}, the type of a consumer). 

We denote by $Y$, $Z_1$, $Z_{-1}$, $X$, $\epsilon_1$, $\epsilon_{-1}$, $\gamma$, and  $\beta$ the population quantities corresponding to the lower cases letters indexed by $i$. The random vectors $X$ and $\beta$ are elements of the unit sphere $\xSd$ of $\xRd$. For the main results of this paper we maintain the following restrictions on the distribution of $(\beta^\top,X^\top)^\top$. 
\begin{assumption}\label{ass1}
\begin{enumerate}[\textup{(}{A1.}1\textup{)}] 
\item\label{ass11} $X$ and $\beta$ are independent,
\item\label{ass12} $X$ and $\beta$ have densities $f_X$ and $\fbeta$ with respect to the spherical measure $\sigma$.
\end{enumerate}
\end{assumption}
\begin{assumption}\label{ass2}
\begin{enumerate}[\textup{(}{A2.}1\textup{)}] 
\item\label{ass13} $\fbeta(x)\fbeta(-x)=0$ for a.e. $x$ in $\xSd$,
\item\label{ass14} The support of $X$, denoted by $\supp(f_X)$, is $H^+=\{x\in\xSd:\
x_1\ge0\}$,
\item\label{ass15} $f_X$ is known and we have $A_X\myeq\|f_X\|_{\xLinfty(H^+)}<\infty$ and $B_X\myeq\|1/f_X\|_{\xLinfty(H^+)}<\infty$.
\end{enumerate}
\end{assumption}
Under Assumption \ref{ass1}, $\fbeta$ is solution of the ill-posed inverse problem: for a.e. $x\in H^+$
\begin{equation}\label{e20}
\E[Y|X=x]=\int_{\xSd}{\rm sign}\left(x^\top y \right)f_{\beta}(y)d\sigma(y)\myeq\mathcal{K}\fbeta(x).
\end{equation}
The operator $\mathcal{K}$ in \eqref{e20} is a convolution on $\xSd$. Estimation of $\fbeta$ in \eqref{e20} is thus related to statistical deconvolution on $\xSd$ (see, \emph{e.g.}, \cite{HK,KK,KPNP}). 
However, the left-hand side of \eqref{e20} is not a density but a regression function where the regressors are random. The identification issue in this model stems from the fact that:
(1) the distribution of the observed data only characterizes $\mathcal{K}\fbeta$ on $\supp(f_X)$ which is a proper subset of $\xSd$ and (2) due to the ${\rm sign}$ function $\mathcal{K}$ has an infinite dimensional null space.  The support of $X$ can only be as large as $H^+$ because the first coordinate of $X$ is positive. This is because we allow for the term $\epsilon_{1,i}-\epsilon_{-1,i}$ in \eqref{eqnet}.

A simple estimator for the density of $\beta$ in this model is given in \cite{GK}. There, rates of convergence for the $\xLn^p$-losses for $1\le p\le \infty$ over Sobolev ellipso\"ids based on the same $\xLn^p$ space (as well as confidence intervals for the value of the density at a point, treatment of endogenous regressors, and of models where some coefficients are nonrandom) are obtained under similar assumptions for choices of the smoothing parameters which depend on unknown parameters of the Sobolev ellipso\"ids. It is assumed in \cite{GK} that the support of $\beta$ lies in an (unknown) hemisphere, namely, that there exists $n$ (unknown) in $\xSd$ such that $\mathbb{P}(n^\top\beta>0)=1$.  This assumption first appeared in \cite{IT} and is stronger than (A2.\ref{ass13}). It implies that for some difference of the characteristics, or taking a limit of these, everyone chooses the same alternative. 
In contrast, (A2.\ref{ass13}) is much less restrictive and does not imply "unselected samples". 
However, everything in \cite{GK} also holds under (A2.\ref{ass13}). 
Assumption (A2.\ref{ass14}) requires that the support of $Z_{1}-Z_{-1}$ is $\xR^{d}$ and is also made in \cite{GK,IT}. 
\cite{GH} allows for continuous regressors which support is a proper subset at the expense of assuming some form of 
unselected samples and relying on integrability assumptions involving $f_{\beta}$It is possible to obtain identification of $f_{\beta}$ when we relax (A2.\ref{ass14}) and the requirement that $f_X$ exists (see (A1.\ref{ass12})). 
This is done in \cite{GG}. The estimation in this case is the subject of future work.
(A2.\ref{ass15}) strengthens (A2.\ref{ass14}) and is used to obtain rates of convergence. It could be viewed as an assumption on the tails of $X$. It is relaxed in \cite{GK} and in this paper at the end of Section \ref{s4}.  
Note as well that Assumption (A1.2) allows for one nonrandom coefficient in the original scale and that  when there are more than two, one should proceed as in Section 5.2 in \cite{GK} with the estimator developed in this paper.

In this paper, we show that the estimator in \cite{GK} can be written as a plug-in of a linear needlet estimator. 
Needlets are a class of linear combinations of spherical harmonics which 
form a tight frame of localized functions on spheres (see \cite{NPW1}). 
Hard-thresholding of series estimators based on needlets have been successfully used in statistics 
for estimation of functions defined on spheres
(see \cite{BKMP} for densities, \cite{Monnier} for regression functions, and  \cite{KKPPP,KPPW,KPNP} for some inverse problems) or compact manifolds (see \cite{KNP}).
This paper proves lower bounds on the minimax risk when the degree of integrability in the loss - specified by the statistician - can differ from the degree of integrability of the Besov body containing the unknown $\fbeta$, giving rise to sparse and dense regimes. The lower bounds correspond, up to logarithmic factors, to the upper bounds in \cite{GK} over Sobolev ellipso\"ids and matching degrees of integrability. This paper proposes to replace the linear needlet estimator in \cite{GK} by a nonlinear estimator based on hard-thresholding with data-driven thresholds and use the same plug-in strategy as in \cite{GK}. The upper bounds on the risk of the estimator also correspond to the lower bounds up to a logarithmic factor, but over all Besov bodies, including nonmatching degrees of integrability. Both the upper and lower bounds are also given for the sup-norm loss.
The data-driven thresholds are similar in spirit to \cite{BLPR} for
density estimation using the Dantzig selector (see also
\cite{DGM,Monnier} for other local thresholding procedures over the
sphere), they are based on sharp concentration inequalities and make
the implementation of the estimator feasible as it is independent of
features of the unknown density.
Proofs are given in the appendix.

\section{Preliminaries}\label{s2}
We use the notation $x\wedge y$ and $x\vee y$ for the minimum and the maximum between $x$ and $y$.  We write
$x\lesssim y$ when there exists $c$ such that $x\le c y$,
$x\gtrsim y$ when there exists $c$ such that $x\ge c y$, and $x\simeq y$ when $x\lesssim y$ and $x\gtrsim y$. 
We denote by $|A|$  and $\ind_A$ the cardinal and indicator of the set $A$, by $\mathbb{N}$ the nonnegative integers, by $\mathbb{N}^*$ the positive integers, by a.e. almost every, and by a.s. almost surely.
We denote for $1\le p\le\infty$ by $\|\cdot\|_{\ell^p}$ the $\ell^p$-norm of a vector, by $\|\cdot\|_p$ the usual norm
on the space $\xLn^p(\xSd)$ of $p$ integrable real-valued functions with respect to the
spherical measure $\sigma$. We write $\xLn^p_{{\rm odd}}(\xSd)$ (resp. $\xLn^p_{{\rm even}}(\xSd)$) the closure in 
$\xLn^p(\xSd)$ of continuous functions on $\xSd$ which are odd
(\emph{i.e.}, for every $x\in\xSd$, $f(-x)=-f(x)$) (resp. even). Every
$f\in\xLn^p(\xSd)$ can be uniquely decomposed as the sum of an odd and
even function $f^-$ and $f^+$ in $\xLn^p(\xSd)$. The space
$\xLtwo(\xSd)$ is a Hilbert space with the scalar product $\langle\ ,\
\rangle$ derived from the norm, there $f^-$ and $f^+$ are
orthogonal. $\mathcal{D}$ is the set of densities and, as it will become clear after Proposition \ref{propK}, $\nu(d)=d/2$ is the degree of ill-posedness of the inverse problem.

\subsection{Harmonic analysis}
\label{s21}
The basic element 
is the orthogonal decomposition
$
\xLtwo(\xSd)=\bigoplus_{k\in\xN}H^{k,d}
$,
where $H^{k,d}$ are the eigenspaces of the Laplacian $\Delta$
on $\xSd$, corresponding to the eigenvalues $-\zeta_{k,d}$, given by $\zeta_{k,d} \myeq
k(k+d-2)$, of dimension
$
L(k,d)\myeq(2k+d-2)(k+d-2)!/(k!(d-2)!(k+d-2))$.
The space $H^{k,d}$ is spanned by an orthonormal
basis 
$\left(h_{k,l}\right)_{l=1}^{L(k,d)}$ and $H^{0,d}$ by $1$. We also have $\xLn^2_{\odd}(\xSd)=\bigoplus_{p\in\xN}H^{2p+1,d}$ and 
$\xLn^2_{\even}(\xSd)=\bigoplus_{p\in\xN}H^{2p,d}$.  The projector
$L_{k,d}$ onto $H^{k,d}$ is the operator with kernel
\begin{equation}L_{k,d}(x,y)=\sum_{l=1}^{L(k,d)}h_{k,l}(x)h_{k,l}(y)= 
\frac{L(k,d)}{\sigma(\xSd)P_k^{\mu(d)}(1)}P_k^{\mu(d)}\left(x^\top y\right)
\label{ekernel20},
\end{equation}
where 
$\mu(d)=(d-1)/2$, the surface of $\xSd$ is $\sigma(\xSd)=2\pi^{d/2}/\Gamma(d/2)$, and $C_k^{\mu}$ are the Gegenbauer polynomials.
The Gegenbauer polynomials, defined for $\mu>-1/2$, are
orthogonal in the space of square integrable functions on $[-1,1]$ with measure $(1-t^2)^{\mu-1/2}dt$.  We have $P_0^{\mu}(t)=1$, $P_1^{\mu}(t)=2\mu t$
for $\mu\ne0$, $P_1^{0}(t)=2t$, and for every $k\in\mathbb{N}$
\begin{equation}\label{erecursion}
(k+2)P_{k+2}^{\mu}(t)=2(\mu+k+1)tP_{k+1}^{\mu}(t)-(2\mu+k)P_k^{\mu}(t).
\end{equation}
Clearly,
for $f\in\xLn^2(\xSd)$, 
we have $f=\sum_{k=0}^{\infty}L_{k,d}f$ and, due to \eqref{ekernel20}, 
\begin{equation}\label{el2}
\forall x\in\xSd,\ \|L_{k,d}(x,\cdot)\|_2^2=\sum_{l=1}^{L(k,d)}|h_{k,l}(x)|^2
=\frac{L(k,d)}{\sigma(\xSd)}.\end{equation}

Powers $\left(-\Delta\right)^{s}f$ for $s\in\xR$ and $f$ in a Banach space $E_1$ are defined in a Banach space $E_2$ when $L_{k,d}f$ is defined in $E_2$ and 
$\left(-\Delta\right)^{s}f\myeq \sum_{k=0}^{\infty}\zeta_{k,d}^{s}L_{k,d}f$
converges in $E_2$. The best approximation in $\xLn^r(\xSd)$  of a function $f$ by harmonics of degree less or equal to $m$ is
\[
E_{m}(f)_r=\inf_{P\in\bigoplus_{k=0}^mH^{k,d}}\left\|f-P\right\|_r.
\]
\begin{definition}\label{def1}
For $s>0$ and $1\le r\le\infty$, $f$ belongs to the Sobolev space  $\xWn^s_r(\xSd)$ if 
\[
\|f\|_{r,s}=\|f\|_r+\left\|\left(-\Delta\right)^{s/2}f\right\|_r<\infty.
\]
\end{definition}
We denote by $\xWn_{r\ {\rm odd}}^{s}(\xSd)$ the restriction of $\xWn_{r}^{s}(\xSd)$ to odd functions.

\begin{definition}\label{def2}
For $s>0$, $1\le r\le\infty$, and
$0<q\le \infty$, $f$ belongs to the Besov space $B^s_{r,q}(\xSd)$ if
\[
\left\|f\right\|^A_{B^s_{r,q}}=\|f\|_r+\left\|\left(2^{js}E_{2^j}(f)_r\right)_{j\in\xN}\right\|_{\ell^q}
<\infty.
\]
\end{definition}

\subsection{The operator}\label{s21b}
\begin{proposition}\label{propK}
The operator $\mathcal{K}$ satisfies the following properties:\begin{enumerate}[\textup{(}{P1.}1\textup{)}]
\item\label{P10} For every $f\in\xLn^1(\xSd)$, $\mathcal{K}f=\mathcal{K}(f^-)$,
\item\label{P11} If $\mathcal{K}f=\mathcal{K}g$ with $f,g\in\xLn^1_{{\rm odd}}(\xSd)$ then $g=f$,
\item\label{P12}  For every $1\le r\le \infty$, $$\xWn^{\nu(d)+|1/r-1/2|(d-2)}_{r\ {\rm odd}}(\xSd)\subseteq \mathcal{K}(\xLn^r_{{\rm odd}}(\xSd))\subseteq \xWn^{\nu(d)-|1/r-1/2|(d-2)}_{r\ {\rm odd}}(\xSd),$$ where the exponents $\nu(d)\pm|1/r-1/2|(d-2)$ cannot be improved,
\item\label{P13} 
For every $1\le r\le \infty$, there exists $B(d,r)$ such that
\begin{equation}\label{eBernstein2}
\forall K\in\xN,\ \forall P\in\bigoplus_{\substack{ k=0\\k\ {\rm odd}}}^KH^{k,d},\ \|\mathcal{K}^{-1}P\|_{r}\le B(d,r)K^{\nu(d)}\|P\|_{r}.
\end{equation}
\end{enumerate}
Moreover, $\mathcal{K}$ is a self-adjoint and compact operator on $\xLn^2(\xSd)$ with null space 
$\xLn_{{\rm even}}^2(\xSd)$, nonzero eigenvalues $(\lambda_{2p+1,d})_{p\in\xN}$ corresponding to the eigenspaces $H^{2p+1,d}$ for $p\in\mathbb{N}$ 
\begin{equation*}
\lambda_{1,d}=\frac{2|\mathbb{S}^{d-2}|} {d-1},\ \forall p \in \mathbb N^*\
\lambda_{2p+1,d}=\frac{2(-1)^{p}|\mathbb{S}^{d-2}|1\cdot3\cdots(2p-1)}
{(d-1)(d+1)\cdots(d+2p-1)}.
\end{equation*}
For every $d\in\xN\setminus\{1\}$, for every $p\in\xN$, there exists $c_{\lambda}(d),C_{\lambda}(d)>0$ such that 
\begin{equation}\label{eqvp}
c_{\lambda}^{-1}(d) p^{-\nu(d)}\le\left|\lambda_{2p+1,d}\right|\le C_{\lambda}(d) p^{-\nu(d)}.
\end{equation}
$\mathcal{K}$ is a homeomorphism between $\xLn_{{\rm odd}}^2(\xSd)$ and 
$\xWn_{2\ {\rm odd}}^{\nu(d)}(\xSd)$. 
\end{proposition}
The fact that $\nu(d)$ is the degree of ill-posedness of the inverse problem follows from (P1.\ref{P13}) and what follows, in particular \eqref{eqvp}. 

Proposition \ref{propK} implies that every $R\in\xWn_{2\ \odd}^{\nu(d)}(\xSd)$ has a unique inverse given by
\begin{equation}\label{eInv}
\mathcal{K}^{-1}\left(R\right)= \sum_{k\ {\rm odd}}
  \frac{1}{\lambda_{k,d}} L_{k,d}\left(R\right) =
\sum_{k\ {\rm odd}}
\frac{1}{\lambda_{k,d}}
\sum_{l=1}^{L(k,d)}  \langle R, \harm_{k,l}\rangle \harm_{k,l}.
\end{equation}

\subsection{Needlets}\label{s22}
Smoothed projection operators (see \cite{GK}) have good
approximation properties in all $\xLn^p(\xSd)$ spaces and are
uniformly bounded from $\xLn^p(\xSd)$ to $\xLn^p(\xSd)$. 
One such operator, the delayed
means, is the integral operator with kernel
\begin{equation}\label{ekernel}
K^{a,J}(x,y)\myeq\sum_{k=0}^{\infty}a\left(\frac{k}{2^{J}}\right)L_{k,d}(x,y),
\end{equation}
where $J$ is an integer, $a$ is a $\xCn^{\infty}$ and decreasing function on $[0,\infty)$ supported
on $[0,2]$ such that, for every $0\le t\le 2$, $0\le a(t)\le 1$ and, for every
$0\le t\le 1$, $a(t)=1$.  The delayed means operator exhibits nearly exponential localization (see Theorem 2.2 in
\cite{NPW1}) and is a building block for the construction of needlets.

Define $b$ such that
$b^2(t)=a\left(t\right)-a(2t)$ for $t\ge0$.
It is nonzero only when $1/2\le t\le2$, satisfies $b^2(t)+b^2(2t)=1$ for  $1/2\le t\le 1$ and thus
for every $t\ge1$, $\sum_{j=0}^{\infty}b^2\left(\frac{t}{2^j}\right)=1$, also $b^2(t)=a(t)$ for $1\le t\le 2$. Take $a$ such that $b$ is bounded away from 0 on $3/5\le t\le 5/3$. 

The second ingredient for the construction of needlets is a quadrature formula (Corollary 2.9 of \cite{NPW1}) 
with positive weights $\left(\omega(j,\xi)^2\right)_{\xi\in\Xi_j}$ and nodes $\xi\in\Xi_j$ which
integrates functions in
$\bigoplus_{k=0}^{2^{j}}H^{k,d}$ and satisfy, for a constant $C_{\Xi}$ which depends on $d$,
\begin{align*}
\forall j\in\xN,\ \forall\xi\in\Xi_j,\quad&C_{\Xi}^{-1}2^{j(d-1)}\le|\Xi_j|\le
C_{\Xi}2^{j(d-1)}\\
&C_{\Xi}^{-1}2^{-j(d-1)/2}\le \omega(j,\xi)\le C_{\Xi}2^{-j(d-1)/2}.
\end{align*}

Needlets are defined as
\begin{align}
\psi_{j,\xi}(x)&\myeq\omega(j,\xi)\sum_{k=0}^{\infty}
b\left(\frac{k}{2^{j-1}}\right)L_{k,d}(\xi,x)\quad  {\rm if}\ j\in\xN,\ \xi\in\Xi_j,\label{eq:psi}\\
\psi_{0,\xi}(x)&\myeq L_{0,d}(\xi,x)\label{eq:psi0}.
\end{align}
For $j=0$, $\psi_{0,\xi}(x)$ is constant and $\Xi_0$
is a singleton. 

The $\xLn^p$-norms of the needlets satisfy, for a constant $\ConstNeedNorm_p$ that can depend on $d$,
\begin{equation}\label{eq:CtrlNormNeed}
\forall j\in\xN,\ \forall\xi\in\Xi_j,\ \ConstNeedNorm_p^{-1} 2^{j(d-1)(1/2-1/p)} \leq \|\psi_{j,\xi}\|_p \leq \ConstNeedNorm_p 2^{j(d-1)(1/2-1/p)}.
\end{equation}
If $f\in\xLn^p(\xSd)$ for $1\le p\le \infty$, then
$
f=\sum_{j=0}^{\infty}\sum_{\xi\in\Xi_j} \langle f, \psi_{j,\xi}\rangle
\psi_{j,\xi}
$.
The needlets form a tight frame, with unitary tightness constant, this means that for $f\in\xL^2(\xSd)$ \[
\|f\|_2^2=\sum_{j=0}^{\infty}\sum_{\xi\in\Xi_j} \left|\langle f,
  \psi_{j,\xi}\rangle\right|^2.
\]
Needlets do not form a basis and there is redundancy. 
Lemma \ref{lem1} (see \cite{BKMP}) relates $\xLn^p(\xSd)$ norms at level $j$ to $\ell^p$ norms of needlet coefficients. Constants may depend on $d$.
\begin{lemma}\label{lem1} 
\begin{enumerate}[\textup{(}i\textup{)}]
\item\label{lem1i} For every $1\le p\le \infty$, there exists a constant $\ConstNeedEquiv_p$ such that for every $j\in\xN$ and $(\beta_{\xi})_{\xi\in\Xi_j}\in \xR^{\Xi_j}$ 
\begin{equation}\label{eNeedEquiv1}
\left\| \sum_{\xi\in\Quad_j} \beta_{\xi} \psi_{j,\xi}\right\|_p \le  \ConstNeedEquiv_p 2^{j(d-1)(1/2-1/p)} \left\|\left(\beta_{\xi}\right)_{\xi\in\Xi_j}\right\|_{\ell^p},
\end{equation}
\item\label{lem1b} There exists constants $c_A$ and $c_{p,A}$ and sets $A_j\subset\Xi_j$
with $|A_j|\ge c_A2^{j(d-1)}$ for $j\in\xN$ such that for every
$1\le p\le \infty$, $j\in\xN$, and $(\beta_{\xi})_{\xi\in A_j}\in\xR^{A_j}$,
\begin{equation}\label{eNeedEquiv1b}
\left\| \sum_{\xi\in A_j} \beta_{\xi} \psi_{j,\xi}\right\|_p \ge  c_{p,A} 2^{j(d-1)(1/2-1/p)} \left\|\left(\beta_{\xi}\right)_{\xi\in A_j}\right\|_{\ell^p},
\end{equation}
\item\label{lem1ii} For every $1\le p\le \infty$, there exists a constant $\ConstNeedEquivTwo_p$ such that for every $j\in\xN$
\begin{equation}\label{eNeedEquiv2}\left( \sum_{\xi\in\Quad_j} | \langle f,
   \psi_{j,\xi}\rangle|^p \right)^{1/p}2^{j(d-1)(1/2-1/p)}
\leq \ConstNeedEquivTwo_p \|f\|_p.
\end{equation}
\end{enumerate}
\end{lemma}

Needlets are such that (see \cite{NPW1}), for all function $a$ in the definition of the smoothed projection operators, the norm
$\left\|\cdot\right\|^A_{B^s_{r,q}}$ defining the Besov spaces is equivalent
to 
\[
\left\|f\right\|_{B^s_{r,q}}=\left\|\left(2^{j(s+(d-1)(1/2-1/r))}\left\|
      \left(\langle f,
        \psi_{j,\xi}\rangle\right)_{\xi\in\Xi_j}\right\|_{\ell^r}\right)_{j\in\xN}\right\|_{\ell^q}.
\]
The ball of radius $M$ for this norm is denoted by $B^s_{r,q}(M)$. 

Recall the following consequence of the proof of the continuous embeddings in
\cite{BKMP}.
\begin{lemma}\label{lem2}
\begin{enumerate}[\textup{(}i\textup{)}]
\item\label{lem2i} If $p\le r\le\infty$, then we have $B^s_{r,q}(M)\subseteq B^s_{p,q}(C_{\Xi}^{1/p-1/r}M)$,
\item\label{lem2ii} If $s>(d-1)(1/r-1/p)$ and $r\le p\le\infty$, then we have $B^s_{r,q}(M)\subseteq B^{s-(d-1)(1/r-1/p)}_{p,q}(M)$,
\item\label{lem2iii} If $f\in B^s_{r,q}(M)$ and $\left(\beta_{j,\xi}\right)_{\xi\in\Xi_j, j\in\xN}$ are its needlet coefficients, then there exists $(D_{j})_{j\in\xN}\in\xR^{\xN}$ such that $\|(D_{j})_{j\in\xN} \|_{\ell^q}\le M$ and
\begin{equation}\label{eq:CtrlScale}
\forall z\ge1,\  \forall j\in\xN,\ \sum_{\xi\in\Quad_j}
|\beta_{j,\xi}|^z
\leq \ConstQuadNb^{1-\infzr/r} D_{j}^{z}
2^{-jz(s+(d-1)(1/2-1/\infzr))}.
\end{equation}
\end{enumerate}
\end{lemma}
Finally recall that, when $f\in B^s_{r,q}$ with
$s>(d-1)/r$, then $f$ is continuous.

\section{Identification of $\fbeta$}\label{s24}
Let us present the arguments for the identification of $\fbeta$.
Proposition \ref{propK} (P1.\ref{P10}) implies that $\mathcal{K}\fbeta=\mathcal{K}\fbeta^-$ is odd.
Thus under (A2.\ref{ass14}) we can define the odd function $R$ as
\begin{equation}\label{eR}
R(x)=\left\{\begin{array}{ll}
              \E[Y|X=x]&\mbox{for\ a.e.}\ x\in H^+\\
              -\E[Y|X=-x]&\mbox{for\ a.e.}\ x\in -H^+\\
            \end{array}\right.
\end{equation}
and we have, for a.e. $x\in\xSd$, $R(x)=\mathcal{K}\fbeta^-(x)$. Uniqueness of $\fbeta^-$ follows from (P1.\ref{P11}). 
Using, for a.e. $x\in\xSd$ $\fbeta(x)\ge0$ and
$\fbeta^-(x)=(\fbeta(x)-\fbeta(-x))/2$, and condition (A2.\ref{ass13}), yields that, for a.e. $x\in\xSd$, we have
\begin{equation}\label{fund}
\fbeta(x)=2\fbetam(x)\ind_{\fbetam(x)>0}.
\end{equation}

In this paper we normalize the vectors of random coefficients and covariates to have unit norm. Indeed, since only the sign of the net utility \eqref{eqnet} matters for choosing between 1 and -1 and the index is linear, a scale normalization of $(\epsilon_{1}-\epsilon_{-1}, \gamma^\top)$ is in order. 
Let us compare with the normalization in \cite{GH}. It is based on the following assumption, which is stronger than the condition in \cite{IT}, that the support of $\beta$ is a subset of some (unknown) hemisphere, which itself is stronger than (A2.\ref{ass13}).
\begin{center}
(H): a.s. there exists $j\in\{1,\dots,d\}$, the coordinate $\gamma_j$ of $ \gamma$ has a sign (excluding 0).
\end{center}
Assumption (H) is likely to hold when $Z_{1j}$ and $Z_{-1j}$ are cost factors, since consumers dislike an increase in cost. 
If (H) holds we can identify for which index $j$ $\gamma_j$ has a sign since it amounts to the finding for which coordinate $z_j$ of $z$ 
$z_j\to \E[Y|Z_1-Z_{-1}=z]$ is (globally) monotone. We can identify the sign of the coefficient by assessing whether the function is increasing (positive) or decreasing (negative).
If $\gamma_j>0$ then we normalize the vector of coefficients by dividing by $\gamma_j$.
If $\gamma_j<0$ we change the sign of $Z_{1j}-Z_{-1j}$ to make it positive. 
A potential issue with this normalization is that if $\beta_j$ can take small values then estimators could differ in finite samples depending on which coefficient is used for normalization.  
Also, monotonicity in one regressor of the conditional mean function implies a type of weak monotonicity (in the sense used to identify treatment effects, see, \emph{e.g.}, \cite{GH}) at the individual level as we now explain. 
Assuming that $\gamma_j>0$, $z_{1i}-z_{-1i}=z$ for all $i=1,\dots,n$, and that we change $z_j$ to $z_j'>z_j$ while leaving unchanged $(\epsilon_{1i}-\epsilon_{-1i},\gamma_i^\top)$ (the characteristics of the individuals) and the other components of $z$, then some people do not change their decision and some choose alternative 1 while originally they had chosen alternative -1, but no one changes from alternative 1 to alternative -1. Monotonicity of the conditional mean function implies monotonicity for every individual. 
This is sometimes not a realistic model of individuals making choices. Clearly (A2.\ref{ass13})  allows both individuals to switch from 1 to -1 and individuals to switch from -1 to 1 after similar changes in $z$ (or $x$). 
On the other hand, if (H) holds then (A2.\ref{ass14}) can be relaxed and we can consider an index which is nonlinear in $X$ (\emph{cf.} \cite{GH}).

\section{Lower bounds}\label{s3}
We take $1\le p,r\le \infty$, $0\le q\le\infty$, $z\ge 1$, and $s>0$,  and consider the minimax risk
\begin{equation}\label{eminimax}\mathcal{R}_n^*\myeq\inf_{\hatfbeta}\sup_{\fbeta\in B^s_{r,q}(M)\cap\mathcal{D}}
\Esp\left\|\hatfbeta-\fbeta\right\|_p^z,
\end{equation}
where the infimum is over all estimators based on the i.i.d. sample of size $n$. The degree of integrability $r$ in the smoothness class $B^s_{r,q}(M)$ is allowed to differ from the degree of integrability $p$
in the loss function. We distinguish two zones for $s,r,q,d$, and $p$:\\ 
(1) the \emph{dense zone} where $s\ge p\left(\nu(d)+(d-1)/2\right)\left(1/r-1/p\right)$ with the restriction 
$q\le r$ if $s=p\left(\nu(d)+(d-1)/2\right)\left(1/r-1/p\right)$, where the rate involves
$$\mu_{{\rm dense}}(d,p,r,s)\myeq s/(s+\nu(d)+(d-1)/2),$$
(2) the \emph{sparse zone} where $(d-1)/r<s< p\left(\nu(d)+(d-1)/2\right)\left(1/r-1/p\right)$, where the rate involves
$$\mu_{{\rm sparse}}(d,p,r,s)\myeq (s-(d-1)(1/r-1/p))/(s+\nu(d)-(d-1)(1/r-1/2)).$$
The terminology dense and sparse is justified by the following heuristic.
The proofs of the lower bounds replace the infimum in \eqref{eminimax} by a minimum over a set of functions which are difficult to estimate.
The functions used to prove the lower bound in the dense zone are functions which could have many nonzero needlet coefficients for $\xi\in A_j$ (see Lemma \ref{lem1}) and a well-chosen $j$. Those used to prove the lower bound in the sparse zone only have two nonzeros. In the dense zone, the rate is the same as for the matched case when $r=p$ studied in \cite{GK}. 
\begin{theorem}\label{t1}
\begin{enumerate}[\textup{(}i\textup{)}]
\item\label{t1i} In the dense zone we have
\begin{equation}\label{eLB1}
\mathcal{R}_n^*\ge c_{{\rm dense}}(d,M,p,r,s,z) \left(\frac{1}{\sqrt{nA_X}}\right)^{\mu_{{\rm dense}}(d,p,r,s)z},
\end{equation}
\item\label{t1ii} In the sparse zone we have
\begin{equation}\label{eLB2}
\mathcal{R}_n^*\ge c_{{\rm sparse}}(d,M,p,r,s,z) \left(\sqrt{\frac{\ln(nA_X)}{nA_X}}\right)^{\mu_{{\rm sparse}}(d,p,r,s)z},
\end{equation}
\end{enumerate}
where the constants $c_{{\rm dense}}$ and $c_{{\rm sparse}}$ depend on $d$, $M$, $p$, $r$, $s$ and $z$.

\end{theorem}
The values of $\mu_{{\rm dense}}$ and $\mu_{{\rm sparse}}$ depend on $d$ through the dimension of $\xSd$.
This is the usual curse of dimensionality in nonparametric regression or density estimation. 
They also depend on $d$ through the degree of ill-posedness $\nu(d)=d/2$ of the inverse problem. 

\section{Adaptive estimation by needlet thresholding}\label{s4}
Consider the estimator 
$
\hatfbeta =2 \hatfbetam
\ind_{\hatfbetam>0}
$, 
where $ \hatfbetam$ is an estimator of $\fbetam$.

\subsection{Smoothed projections and linear needlet estimators}\label{s41}
A smoothed projection estimator of $\fbeta^-$ with kernel \eqref{ekernel}, window $a$, and $J\in\xN$, is given for $x\in\xSd$ by
\begin{align*}
  \hatfbetam^{a,J} (x)= \sum_{k\ {\rm odd}}
  \frac{a\left(\frac{k}{2^{J}}\right)}{\lambda_{k,d}}
    \widehat{L_{k,d}R}(x),
\end{align*}
with the unbiased estimator of $L_{k,d}R(x)$ (see Lemma \ref{lemrandomdesign}): $\widehat{L_{k,d}R}(x)=0$ if $k$ is even, else
\[
\widehat{L_{k,d}R}(x)=\frac{2}{n} \sum_{i=1}^{n} \frac{ y_i
  L_{k,d}(x_i,x)}{\fX(x_i)}.
\]
Alternatively, we can estimate $\fbetam$ using the needlet frame
with smoothing window $a$.
The coefficients $\beta^{a}_{j,\xi}= \langle \fbetam, \psi_{j,\xi} \rangle$ are such that
\begin{align*}
\beta^{a}_{j,\xi}&= \omega(j,\xi) \sum_{k \odd} b\left(\frac{k}{2^{j-1}}\right) \langle \fbetam, L_{k,d}(\xi,\cdot)\rangle\\
&= \omega(j,\xi) \sum_{k \odd}
\frac{b\left(\frac{k}{2^{j-1}}\right)}{\lambda_{k,d}} \langle L_{k,d}R, L_{k,d}(\xi,\cdot)\rangle\\
&= \omega(j,\xi) \sum_{\substack{k \odd\\2^{j-2}<k<2^j}}
\frac{b\left(\frac{k}{2^{j-1}}\right)}{\lambda_{k,d}} L_{k,d}R(\xi).
\end{align*}
Using that $a\left(\frac{k}{2^{j}}\right)=1$ for $k=0,\hdots,2^{j}$ and denoting  by $\fbetam^{a,J}=\E\left[\hatfbetam^{a,J}\right]$,
we obtain that, for $1\le j\le J$,
$\beta^{a}_{j,\xi}=\left\langle \fbetam^{a,J}, \psi_{j,\xi}\right\rangle$, which can be estimated without bias by
\[
\hatbeta^a_{j,\xi}=\omega(j,\xi) \sum_{k \odd}
\frac{b\left(\frac{k}{2^{j-1}}\right)}{\lambda_{k,d}}
\widehat{L_{k,d}R}(\xi)
\stackrel{(\triangle_1)}{=}\left\langle \hatfbetam^{a,J}, \psi_{j,\xi}
\right\rangle.
\]
Moreover, for $x\in\xSd$,
\[
\hatbeta^a_{j,\xi}\psi_{j,\xi}(x)=\omega(j,\xi)^2 \left(\sum_{k \odd}
\frac{b\left(\frac{k}{2^{j-1}}\right)}{\lambda_{k,d}}
\widehat{L_{k,d}R}(\xi)\right)\left(\sum_{k}
b\left(\frac{k}{2^{j-1}}\right)
L_{k,d}(\xi,x)\right)
\]
belongs to $\bigoplus_{k=0}^{2^{j}}H^{k,d}$, thus by the quadrature formula 
\[
\sum_{\xi\in\Xi_j}\hatbeta^a_{j,\xi}\psi_{j,\xi}(x)=\sum_{k \odd}
\frac{b^2\left(\frac{k}{2^{j-1}}\right)}{\lambda_{k,d}}
\widehat{L_{k,d}R}(x).
\]
This yields $\sum_{j=0}^J\sum_{\xi\in\Xi_j}  \hatbeta^a_{j,\xi}  \psi_{j,\xi}=\hatfbetam^{a,J-1}$, indeed
\begin{align*}
\sum_{j=0}^J\sum_{\xi\in\Xi_j}  \hatbeta^a_{j,\xi}  \psi_{j,\xi}
&= \sum_{j=1}^J\sum_{\xi\in\Xi_j}
\hatbeta^a_{j,\xi}  \psi_{j,\xi}\quad\mathrm{(due\ to\ }(\triangle_1)\ \mathrm{and\ because\ } \hatfbetam^{a,J}\ \mathrm{is\ odd)}\\
&\stackrel{(\triangle_2)}{=} \sum_{\substack{1\le k<2^{J-1}\\k \odd}}
\frac{1}{\lambda_{k,d}}
\widehat{L_{k,d}R}+\sum_{\substack{2^{J-1}\le k\le 2^J\\k \odd}}
\frac{b^2\left(\frac{k}{2^{J-1}}\right)}{\lambda_{k,d}}
\widehat{L_{k,d}R}
\\
&\stackrel{(\triangle_3)}{=} \sum_{\substack{1\le k<2^{J-1}\\k \odd}}
\frac{1}{\lambda_{k,d}}
\widehat{L_{k,d}R}
+\sum_{\substack{2^{J-1}\le k\le 2^J\\k \odd}}
\frac{a\left(\frac{k}{2^{J-1}}\right)}{\lambda_{k,d}}
\widehat{L_{k,d}R},
\end{align*}
where
($\triangle_2$) uses that for $1/2\le t\le 1$, $b^2(t)+b^2(2t)=1$, while ($\triangle_3$) that $b^2(t)=a\left(t\right)$ for $1\le t\le 2$. 
Thus, the smoothed projection and needlet estimators coincide. 

\subsection{Nonlinear estimator with data-driven thresholds}\label{s42}Consider, for $\gamma\ge1$ and $\rho_{T_{j,\xi,\gamma}}(x)=x\ind_{|x|>T_{j,\xi,\gamma}}$, the nonlinear estimator of $\fbetam$: 
$$
  \hatfbetam^{a,\rho} =
\sum_{j=0}^J\sum_{\xi\in\Xi_j}  \rho_{T_{j,\xi,\gamma}}\left(\hatbeta^a_{j,\xi}\right) \psi_{j,\xi}.
$$
It is classical that the optimal choice of $J$ for linear estimators depends on the parameters of the smoothness ellipsoid. In contrast, using a 
thresholded estimator allows to take $J$ large and independent of the parameters. 
Thresholding induces additional bias compared to linear estimators which allows to reduce the variance incurred by taking $J$ large.  

The level of thresholding should depend on the size of the coefficients relative to their variance. This variance is proportional to $1/\sqrt{n}$ so that the level of the threshold does not have to depend on the smoothness of the unknown function. Instead of using a conservative upper bound on their variance, as is usually the case in estimation using wavelets, we use data-driven levels of thresholding. These provide better estimators in small samples.
Lemma \ref{toracle} gives a theoretical guarantee that the performance is almost as good as that of an oracle which would know the variance of the estimators of the coefficients. The data-driven thresholding rule uses that
$
\hatbeta^a_{j,\xi}=
 \frac1n\sum_{i=1}^nG_{j,\xi}(x_i,y_i)
$
with
\begin{equation}\label{eqG}
G_{j,\xi}(x_i,y_i)\myeq \frac{2}{n}\sum_{i=1}^{n}
\omega(j,\xi)
\frac{y_i}{\fX(x_i)}
\sum_{k \odd}
\frac{b\left(\frac{k}{2^{j-1}}\right)}{\lambda_{k,d}}
 L_{k,d}(x_i,\xi).
 \end{equation}
Define the estimator of the variance by
\begin{equation}\label{eqvar}
\hat{\sigma}_{j,\xi}\myeq\sqrt{\frac{1}{n(n-1)}\sum_{i=2}^n\sum_{k=1}^{i-1}
\left(G_{j,\xi}(x_i,y_i)-G_{j,\xi}(x_k,y_k)\right)^2},
 \end{equation}
$t_n=\sqrt{\log n/n}$, and the data-driven thresholds
\[
T_{j,\xi,\gamma}\myeq 2 \sqrt{2
\gamma} t_n \widehat{\sigma}_{j,\xi} + \frac{28}{3} M_{j,\xi}
\frac{\gamma \log n}{n-1},
\]
 where $M_{j,\xi}$ is an upper bound
on the sup-norm over $ H^+\times\{\pm1\}$ of $G_{j,\xi}(x,y)-
\Esp \left[
      G_{j,\xi}(X,Y) \right]=G_{j,\xi}(x,y)- \beta^{a}_{j,\xi}$
(\emph{e.g.}, $2
\|G_{j,\xi}\|_{\infty}$). 
For example, using \eqref{eq:CtrlNormNeed} and Proposition \ref{propK}, we
get 
\begin{equation}\label{esup}
2\|G_{j,\xi}\|_{\infty}\le 2
\left\|\mathcal{K}^{-1}\left(\psi_{j,\xi}^-\right)\right\|_{\infty}B_X
\le 2\ConstNeedNorm_{\infty}B(d,\infty)2^{j(\nu(d)+(d-1)/2)}B_X\myeq M_{j}.
\end{equation}
The second term in $T_{j,\xi,\gamma}$ controls the error in estimating the threshold.

\begin{theorem}\label{t2}
For $J$ such that
$2^{J(\nu(d)+(d-1)/2)}B_X^{1/2}\simeq
t_n^{-1}$, $M>0$, and $s>(d-1)/r$, 
\begin{enumerate}[\textup{(}i\textup{)}]
\item\label{t2i}
If $z>1$ and $\gamma>z/2+1$, we have
\begin{equation}
\sup_{\fbeta\in B^s_{r,q}(M)\cap\mathcal{D}}
\Esp\left\|\hatfbeta^{a,\rho} -\fbeta\right\|_\infty^z\le
\tilde{c}(d,\infty,r,s,\gamma)(\log n)^{z-1}M^r \left(B_Xt_n\right)^{\mu_{{\rm sparse}}(d,\infty,r,s) z}.
\label{eUB1}
\end{equation}
\item\label{t2ii} If $p<\infty$ and
$\gamma>p/2$, we have
\begin{equation}
\sup_{\fbeta\in B^s_{r,q}(M)\cap\mathcal{D}}
\Esp\left\|\hatfbeta^{a,\rho} -\fbeta\right\|_p^p\le
\tilde{c}(d,p,r,s,\gamma)(\log n)^{p-1}M^{\varpi} \left(B_Xt_n\right)^{\mu(d,p,r,s)  p},\label{eUB2}
\end{equation}
where $\mu(d,p,r,s) =\mu_{{\rm dense}}(d,p,r,s) $ and $\varpi=r$ in the dense zone,
while $\mu(d,p,r,s) =\mu_{{\rm sparse}}(d,p,r,s) $ and
$\varpi>p\frac{\nu(d)+(d-1)(1/2-1/p)}{s+\nu(d)-(d-1)(1/r-1/2)}$ is arbitrary in the sparse zone, and 
$\tilde{c}(d,p,r,s,\gamma)$ is a constant which depends on $d,p,r,s$, and $\gamma$.
\end{enumerate}
\end{theorem}

The upper bounds in Theorem \ref{t2} match the lower bound in Theorem \ref{t1} up to logarithmic factors. Hence, the proposed estimator is minimax adaptive (up to the $\log$ factors).

\section{Simulation study}
We study the performance of the estimator when $d=3$, 
$n=3000,5000,10000$, and $X$ is uniform on $H^+$.
We use of the Von Mises-Fisher distribution $\mathrm{vMF}(\mu,\kappa)$ with density \begin{equation*}
f(\beta;\mu,\kappa) = \frac{\kappa}{4\pi\sinh\kappa} \exp\left(\kappa\mu^\top \beta\right)
\end{equation*}
with respect to $\sigma$. We take $\beta=(\tilde{\beta_1},\tilde{\beta_2},|\tilde{\beta_3}|)$ in the cases:
\begin{itemize}
	\item $\tilde \beta$ follows a $\mathrm{vMF}(\mu,\kappa)$ distribution 
	where $\mu=(0\ 0\ 1)^\top$
	and $\kappa = 10$.
	\item $\tilde \beta$ follows a mixture 
$\lambda \mathrm{vMF}(\mu_1,\kappa) + (1-\lambda)\mathrm{vMF}(\mu_2,\kappa),
$ where $\mu_1=(2^{-1/2}\ 0\ 2^{-1/2})^\top$, $\mu_2=(-2^{-1/2}\ 0\ 2^{-1/2})^\top$, $\kappa=10$ and $\lambda = 0.3$.
\end{itemize}
We use the cubature defined in spherical coordinates as a product of
the Gauss-Legendre quadrature with $m$ nodes and trapezoid rule with
$2m$ subdivisions (see  \cite{AH}). The resulting cubature has $2m^2$
nodes and integrates exactly all polynomials on the sphere up to
degree $2m-1$. We take the same function $a$ as in \cite{BKMP}.

The threshold is driven by the parameter $\gamma$. The choice of
$\gamma$ slightly depends on the targeted norm. Here we focus on a
simultaneous control of the 
$\xL^1$, $xL^2$, $xL^4$ and $\xL^{\infty}$ norm. According to our
analysis, $\gamma$ should be chosen stricly larger than $4$. We have
nevertheless chosen to use $\gamma=4$ which turns out to be sufficient
in practice.

\begin{figure}
\centering
	\begin{subfigure}[b]{0.35\textwidth}
		\includegraphics[width=\textwidth]{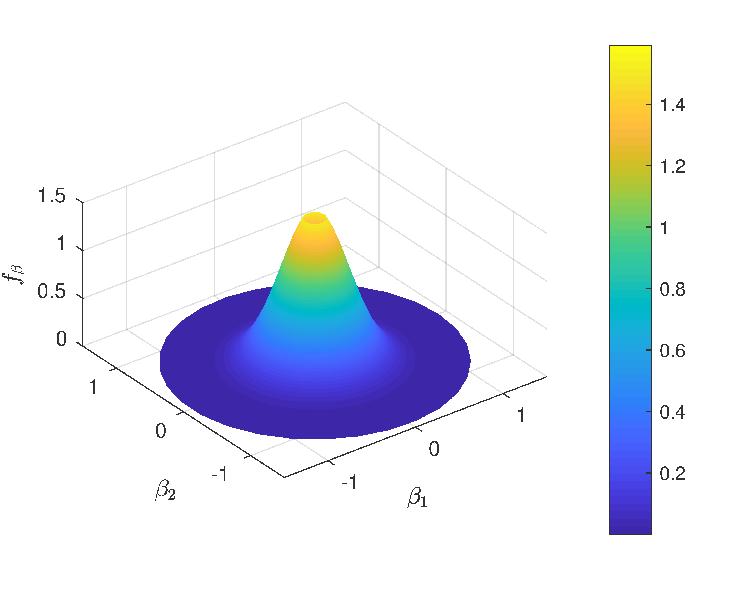}
		\caption{True density}
	\end{subfigure}
	\begin{subfigure}[b]{0.35\textwidth}
		\includegraphics[width=\textwidth]{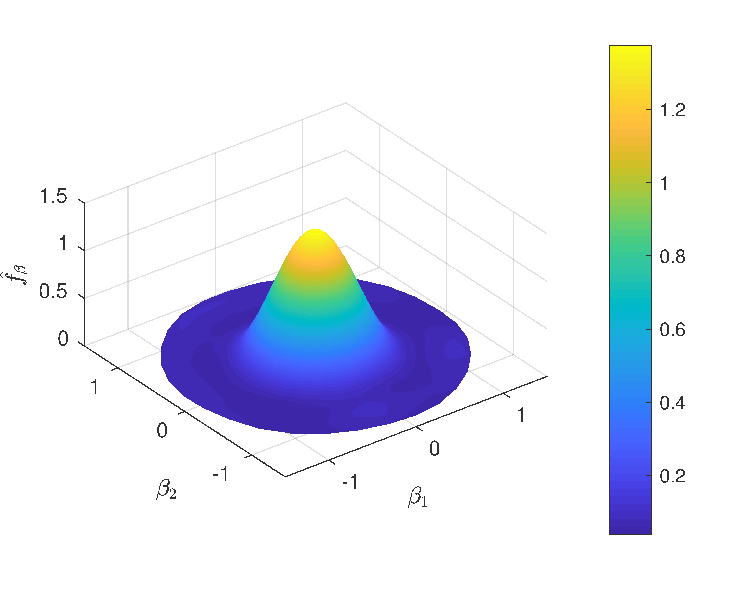}
		\caption{Mean of estimates}
	\end{subfigure}
	\begin{subfigure}[b]{0.35\textwidth}
		\includegraphics[width=\textwidth]{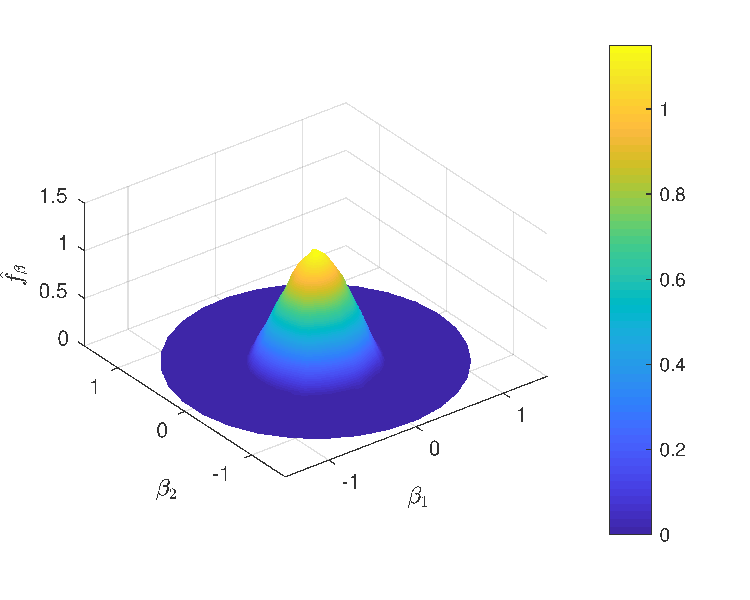}
		\caption{5\% quantile of estimates}
	\end{subfigure}
	\begin{subfigure}[b]{0.35\textwidth}
		\includegraphics[width=\textwidth]{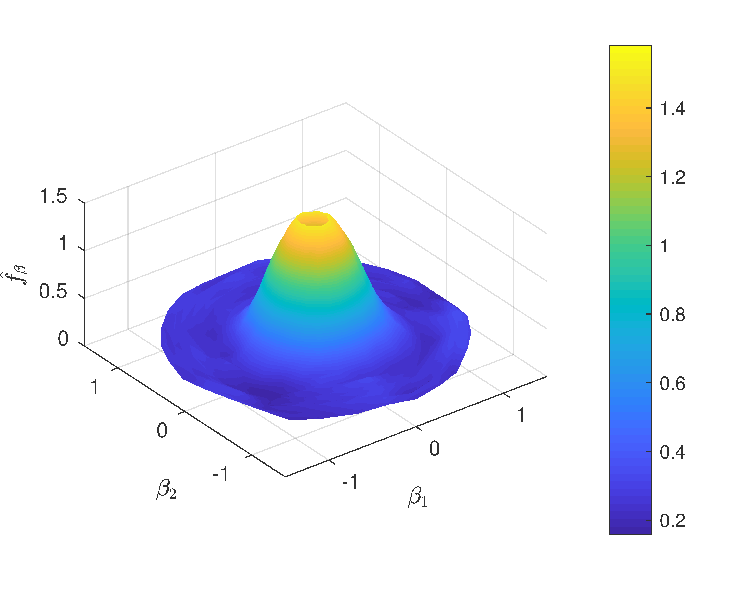}
		\caption{95\% quantile of estimates}
	\end{subfigure}
	\begin{subfigure}[b]{0.35\textwidth}
		\includegraphics[width=\textwidth]{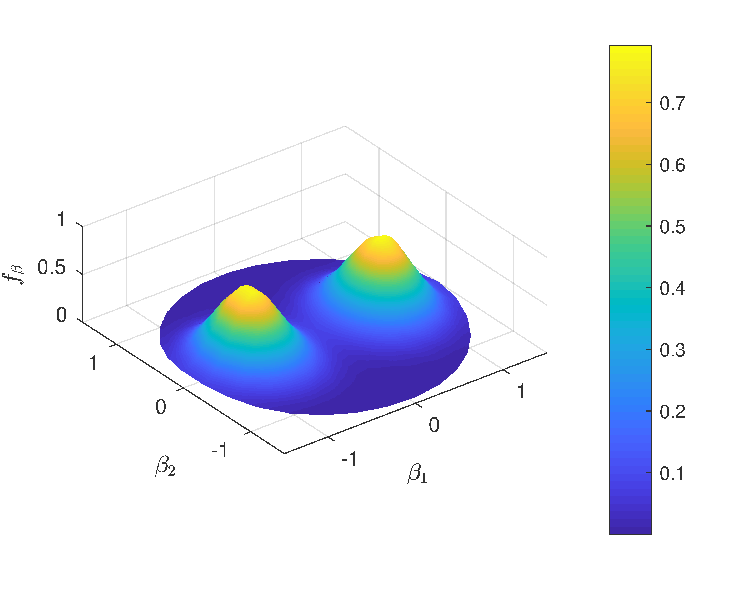}
		\caption{True density}
	\end{subfigure}
	\begin{subfigure}[b]{0.35\textwidth}
		\includegraphics[width=\textwidth]{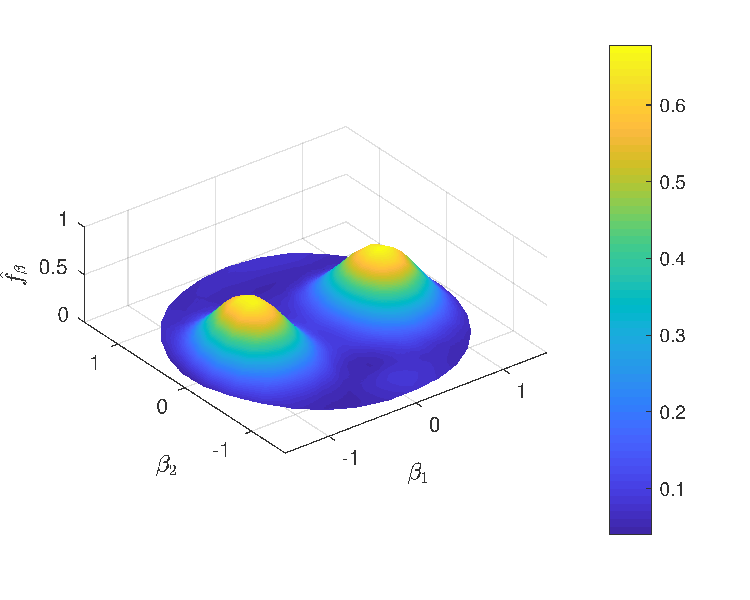}
		\caption{Mean of estimates}
	\end{subfigure}
	\begin{subfigure}[b]{0.35\textwidth}
		\includegraphics[width=\textwidth]{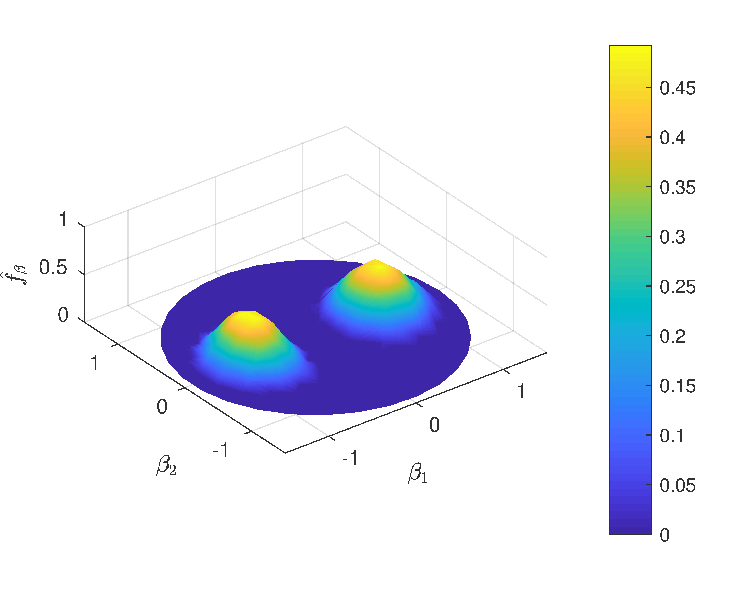}
		\caption{5\% quantile of estimates}
	\end{subfigure}
	\begin{subfigure}[b]{0.35\textwidth}
		\includegraphics[width=\textwidth]{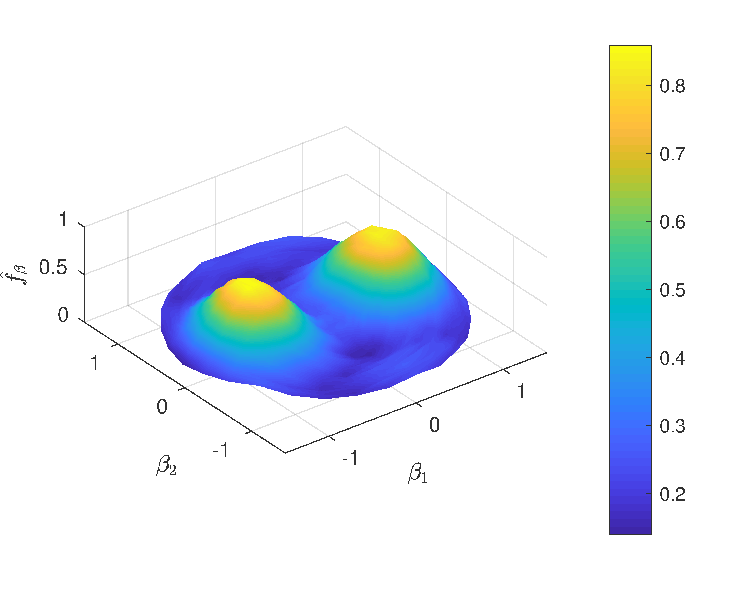}
		\caption{95\% quantile of estimates}
	\end{subfigure}
	\caption{{\small True density and distribution of the estimates.}}
	\label{fig}	
\end{figure}

Figure~\ref{fig} displays the distribution of estimates based on a Monte-Carlo experiments with 100 replications and $n=3000$.  We plot the Lambert equal-area projection on the disk which is defined (see \cite{MJ})
\begin{equation*}
  (\sin\theta\cos\phi,\sin\theta\sin\phi,\cos\theta)^\top \mapsto 2\sin\left(\frac{\theta}{2}\right)(\cos\phi,\sin\phi)^\top.
\end{equation*}

\begin{figure}
  \centering
  \includegraphics[width=12cm]{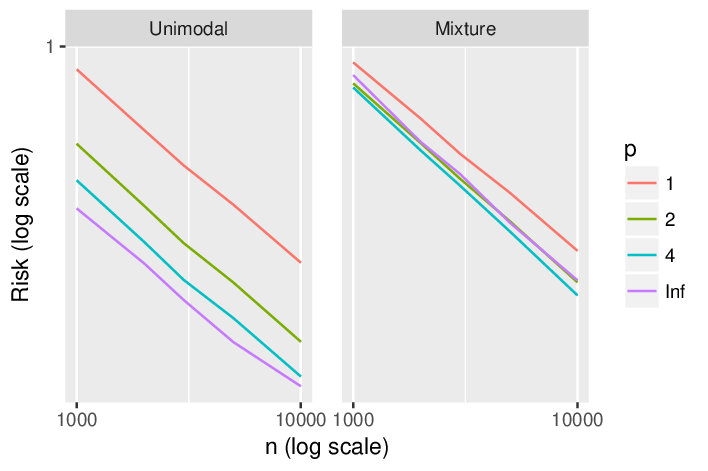}
  \caption{\small Decay of the risk with $n$ in logarithmic scales.}
  \label{riskdecay}
\end{figure}

Our main contribution is a control of the estimation error for
all $\xL^p$ norm. Table~\ref{tab} displays the expected risk, approximated using Monte-Carlo and
100 replications, for some $\xL^p$ norms. More precisely, we have approximated the following renormalized quantities:
$\left(\mathbb{E}\left[\left\|\widehat{f}_{\beta}-f_{\beta}\right\|_p^p\right]/\|f_{\beta}\|_p^p\right)^{1/p}$
for $p = \{1, 2 ,4\}$ and
$\mathbb{E}\left[\left\|\widehat{f}_{\beta}-f_{\beta}\right\|_{\infty}\right]/\|f_{\beta}\|_{\infty}$. Figure~\ref{riskdecay}
displays the decay of those error with respect to $n$ in a logarithmic
scales. As expected, we observe a simultaneous control over all norm
and the error decays follows the power law given by the upper
bounds. The results are similar to the one obtained in \cite{GK}
except that our threshold does not depend on the unknown regularity of
the function whereas the level used in \cite{GK} depends on it.

\begin{table}[H]
\centering
				\begin{tabular}{c|ccccc||ccccc}
						& \multicolumn{5}{c}{Unimodal} & \multicolumn{5}{c}{Mixture} \\\hline
                                  \backslashbox{Risk}{$n$}
                                  \footnotesize	& $1000$ & $2000$ &
                                                                    $3000$ & $5000$ & $10000$ & $1000$ & $2000$ & $3000$ & $5000$ & $10000$ \\
                                  \hline
						{\footnotesize$\mathbb{E}\left[\left\|\widehat{f}_{\beta}-f_{\beta}\right\|_1\right]/\|f_{\beta}\|_1$} &
$0.89$&$0.64$&$0.53$&$0.43$&$0.32$&$0.92$&$0.68$&$0.57$&$0.46$&$0.34$\\
                                  {\footnotesize$\left(\mathbb{E}\left[\left\|\widehat{f}_{\beta}-f_{\beta}\right\|^2_2\right]/\|f_{\beta}\|_2^2\right)^{1/2}$}
                                                &
          $0.6$&$0.43$&$0.35$&$0.29$&$0.21$&$0.82$1&$0.6$&$0.5$&$0.4$&$0.29$\\                                        
			{\footnotesize$\left(\mathbb{E}\left[\left\|\widehat{f}_{\beta}-f_{\beta}\right\|^4_4\right]/\|f_{\beta}\|_4^4\right)^{1/4}$}
                                                &
$0.49$&$0.36$&$0.29$&$0.24$&$0.17$&$0.8$&$0.58$&$0.48$&$0.38$&$0.27$\\
                                                  {\footnotesize$\mathbb{E}\left[\left\|\widehat{f}_{\beta}-f_{\beta}\right\|_{\infty}\right]/\|f_{\beta}\|_{\infty}$}
                                                                                                                         &
$0.4$2&$0.32$&$0.26$&$0.21$&$0.17$&$0.86$&$0.6$&$0.51$&$0.39$&$0.29$
					\end{tabular}
	       \caption{{\small Risk.}}
	       \label{tab}
\end{table}

\section{Appendix}\label{sproofs}
\subsection{A preliminary lemma}
\begin{lemma}\label{lemrandomdesign}
The following equality holds for every $g\in\xLn^2(\xSd)$,
$$\langle R, g \rangle= 2 \Esp \left[ \frac{ Y g^{-}(X)}{\fX(X)}\right].$$
\end{lemma}
\noindent{\bf Proof.} The result is based on the following 
\begin{align*}
  \langle R, g \rangle
&= \langle R, g^- \rangle\quad \mbox{(because $R$ is odd)} \\
&= 2 \int_{H^+}  \frac{R(x)g^{-}(x)}{\fX(x)} \fX(x)
  d\sigma(x)\\
 &= 2 \Esp \left[
\frac{R(X)g^{-}(X)}{\fX(X)}
\right]\\
& = 2 \Esp \left[
\frac{\Esp[Y|X]g^{-}(X)}{\fX(X)}\right].
\qquad\qquad\qquad\qquad\qquad\square
\end{align*}

\subsection{Proof of Proposition \ref{propK}}
The operator $\mathcal{K}$ is related to the Hemispherical transform (see \cite{GK,Rubin}) defined for $f\in\xLn^1(\xSd)$ and a.e. $x\in\xSd$ by
$$\mathcal{H}(f)(x)\myeq\int_{\xSd}\ind_{x^\top y 
>0}f(y)d\sigma(y),$$ through
$$Kf=2\mathcal{H}(f)-\int_{\xSd}f(y)d\sigma(y).$$
(P1.\ref{P10}) is a consequence of the fact that $y\to x^\top y\in\xLn^{\infty}_{{\rm odd}}(\xSd)$. (P1.\ref{P11}) follows from Theorem 2 (ii), and (P1.\ref{P12}) follows from Theorem C in \cite{Rubin}.
The second part of the proposition together with (P1.\ref{P13}) are consequences of the properties of $\mathcal{H}$ detailed in \cite{GK}.  The inequalities \eqref{eqvp} correspond to Lemma A.2. Note 
however that there is a typo in the proof and we should read
$1.3\hdots(2p-1)\asymp p^{-1/2} 2.4\hdots(2p)$ but the result still holds.

\subsection{Proof of Theorem \ref{t1}}\label{s6}
Start by noting that for every $j\in\xN$ and $\xi\in\Xi_j$, 
$$\int_{\xSd}\psi_{j,\xi}(x)dx=\omega(j,\xi)b(0)=\omega(j,\xi)(a(0)-a(0))=0.$$
This implies that the functions $f_m$ that we introduce below integrate to 1.

\subsubsection{Proof of the lower bound in the dense zone}\label{s61}
Consider the family $(P_m)_{m=0}^{\mathcal{M}}$, where $\mathcal{M}\in\xN^*$, of
distributions of an i.i.d. sample of $(Y,X)$ of size $n$ when $\fbeta=f_m$ and the density of $X$ is $\fX$. 
These probabilities are absolutely continuous with respect to the product of $\delta_1+\delta_{-1}$, where $\delta_y$ denotes the Dirac mass at $y$ and $\sigma$. 
Take $j\in\xN$, $f_0=1/\sigma(\xSd)$, and consider the set $A_j$ from Lemma \ref{lem1} \eqref{lem1b}. 
By the Varshamov-Guilbert bound (Lemma 2.9 in \cite{Tsyb}) there exists $\Omega\subseteq \{0,1\}^{A_j}$ containing $(0,\hdots,0)$ such that $|\Omega|=2^{|A_j|/8}$  and $\forall (\omega_1,\omega_2)\in
\Omega^2$, $\|\omega_{1}-\omega_{2}\|_{\ell^1}\ge |A_j|/8$. Enumerate the elements of $\Omega$ from 0 (corresponding to the zero vector) to $\mathcal{M}\myeq |\Omega|-1$ and define  
$$f_m\myeq f_0+\gamma\sum_{\xi\in A_j} \omega_{\xi}\psi_{j,\xi}$$
when $(\omega_{\xi})_{\xi\in A_j}$ is the $m$th element of $\Omega$
and $\gamma=cC_{\Xi}^{-1/r}M2^{-j(s+(d-1)/2)}$ for $0<c<1$ such that all $f_m$ are nonnegative. 
We now use the following result (see Theorem 2.5 in \cite{Tsyb}).
\begin{lemma}\label{lem5}
If for $0<\alpha<1/8$ we have:
\begin{enumerate}[\textup{(}i\textup{)}]
\item\label{lem5i}
$f_m\in B_{r,q}^s(M)\cap\mathcal{D}$ for $m=0,\hdots,\mathcal{M}$,
\item\label{lem5ii} $\forall\ 0\le m<l\le\mathcal{M},\ \|f_m-f_l\|_p\ge 2h>0$,
\item\label{lem5iii} $\frac{1}{\mathcal{M}}\sum_{m=1}^{\mathcal{M}}K(P_m,P_0)\le\alpha \ln(\mathcal{M})$,
\end{enumerate}
then for every $z\ge1$
\begin{equation}\label{eqlb1}
\inf_{\hatfbeta}\sup_{\fbeta\in B^s_{r,q}(M)\cap\mathcal{D}}
\Esp\left\|\hatfbeta-\fbeta\right\|_p^z \ge
h^{z}\frac{\sqrt{\mathcal{M}}}{1+\sqrt{\mathcal{M}}}\left(1-2\alpha-\sqrt{\frac{2\alpha}{\ln(\mathcal{M})}}\right).
\end{equation}
\end{lemma}
Start by checking \eqref{lem5i} in Lemma \ref{lem5}. It is enough to show that $f_m\in B_{r,q}^s(M)$.
Indeed, for $r\ge1$ and $\omega\in\Omega$, we have $\left\|\left(\omega_\xi\right)_{\xi\in
    A_j}\right\|_{\ell^r}\le\left\|\left(\omega_\xi\right)_{\xi\in
    A_j}\right\|_{\ell^1}^{1/r}\le C_{\Xi}^{1/r}2^{j(d-1)/r}$, we obtain
\[
\gamma 2^{j(s+(d-1)(1/2-1/r))}\left\|\left(\omega_\xi\right)_{\xi\in
    A_j}\right\|_{\ell^r}\le \gamma C_{\Xi}^{1/r}2^{j(s+(d-1)/2)} \le M.
\]
Lemma \ref{lem1} \eqref{lem1b} now yields that for every $1\le p\le \infty$ and  $0\le m<l \le \mathcal{M}$
$$\left\|f_{m}-f_{l}\right\|_p\ge \gamma c_{p,A} 2^{j(d-1)(1/2-1/p)} \left(\frac{c_A}{8}2^{j(d-1)}\right)^{1/p}=
2 h.$$
Thus \eqref{lem5ii} in Lemma \ref{lem5} follows with $h= c_{p,A}\left(\frac{c_A}{8}\right)^{1/p}cC_{\Xi}^{-1/r}M 2^{-js-1}$.\\
By independence,
the Kullback-Leibler divergence between $P_m$ and $P_0$ is given by
\[
K(P_m,P_0)=n\Esp\left[\cH(f_m)(X)\ln\left(\frac{\cH(f_m)(X)}{\cH(f_0)(X)}\right)+\left(1-\cH(f_m)(X)\right)
\ln\left(\frac{1-\cH(f_m)(X)}{1-\cH(f_0)(X)}\right)\right].
\]
Using that, for $x>0$, $\ln(x)\le x-1$, we obtain 
\begin{equation*}
K(P_m,P_0)\le n\Esp\left[\frac{\cH(f_m-f_0)(X)^2}{\cH(f_0)(X)\left(1-\cH(f_0)(X)\right)}\right],
\end{equation*}
and thus 
\begin{align*}
K(P_m,P_0)&\le 4nA_X\left\|\cH(f_m-f_0)\right\|_2^2\\
&\le  4nA_X\lambda_{2^j+1,d}^2\left\|f_m-f_0\right\|_2^2,
\end{align*}
where the last display comes from the fact that $f_m-f_0\in\bigoplus_{2^j+1\le k\le 2^{j+2}-1}H^{k,d}$.
From \eqref{eqvp} we get
$$K(P_m,P_0)\le4C_{\lambda}(d)^2nA_X2^{-2j\nu(d)}\left\|f_m-f_0\right\|_2^2,$$
which yields using Lemma \ref{lem1} \eqref{lem1i}
\begin{align*}
K(P_m,P_0)
&\le \left(2C_{\lambda}(d)\ConstNeedEquiv_2\gamma\right)^2nA_X2^{-2j\nu(d)}
\left\|\left(\omega_{\xi}\right)_{\xi\in A_j}\right\|_{\ell^2}^2\\
&\le \left(2C_{\lambda}(d)\ConstNeedEquiv_2\gamma\right)^2nA_X2^{-2j\nu(d)}
\left\|\left(\omega_{\xi}\right)_{\xi\in A_j}\right\|_{\ell^1}\\
&\le \left(2C_{\lambda}(d)\ConstNeedEquiv_2\gamma\right)^2C_{\Xi}nA_X2^{j(d-1-2\nu(d))}\\
&\le \left(2C_{\lambda}(d)\ConstNeedEquiv_2cM\right)^2C_{\Xi}^{1-2/r}nA_X2^{-2j(s+\nu(d))}.
\end{align*}
Condition \eqref{lem5iii} of Lemma \ref{lem5} is satisfied
once
\begin{equation}\label{eqlb2}
\frac{2^5\left(C_{\lambda}(d)\ConstNeedEquiv_2cM\right)^2}{\ln(2)}C_{\Xi}^{-2/r}nA_X2^{-2j(s+\nu(d)+(d-1)/2)}\le\alpha<\frac{1}{8}.
\end{equation}

\noindent For $\alpha<1/8$, the lower bound \eqref{eqlb1} yields that
\begin{align*}
\inf_{\hatfbeta}\sup_{\fbeta\in B^s_{r,q}(M)}
\Esp\left\|\hatfbeta-\fbeta\right\|_p^z &\ge
\left( c_{p,A}\left(\frac{c_A}{8}\right)^{1/p}cC_{\Xi}^{-1/r}M 2^{-js-1}\right)^{z}\left(\frac{3}{4}-\frac{1}{2\sqrt{\ln(\mathcal{M})}}\right)\\
&\ge\frac{1}{2}
\left( c_{p,A}\left(\frac{c_A}{8}\right)^{1/p}cC_{\Xi}^{-1/r}\frac{M}{2}\right)^{z} 2^{-jsz},
\end{align*}
where the inequality leading to the second display holds when $\ln(\mathcal{M})\ge4$, for example 
for $j(d-1)\ge\ln(5/c_A\ln(2))/\ln(2)$. Now \eqref{eqlb2} is satisfied for
$$j\ge j_0\myeq1+\frac{\ln\left(2^8\left(C_{\lambda}(d)\ConstNeedEquiv_2cM\right)^2C_{\Xi}^{-2/r}nA_X/\ln(2)\right)}{2\ln(2)(s+\nu(d)+(d-1)/2)},$$
which implies the lower bound
\begin{align*}
&\inf_{\hatfbeta}\sup_{\fbeta\in B^s_{r,q}(M)}
\Esp\left\|\hatfbeta-\fbeta\right\|_p^z \\
&\ge
\frac{1}{2}
\left( c_{p,A}\left(\frac{c_A}{8}\right)^{1/p}cC_{\Xi}^{-1/r}M2^{-s-1}\right)^{z} 
\left(\frac{2^8\left(C_{\lambda}(d)\ConstNeedEquiv_2cM\right)^2C_{\Xi}^{-2/r}nA_X}{\ln(2)}\right)^{-\mu_{{\rm dense}}(d,p,r,s)z/2}.
\end{align*}

\subsubsection{Proof of the lower bound in the sparse zone}\label{s62}
In this proof we consider asymptotic orders for simplicity. The various constants can be obtained like in Section \ref{s61}. 
Consider the hypotheses
\[
f_m=\frac{1}{\sigma(\xSd)}+\gamma\psi_{j,\xi_m},
\]
where $\xi_m\in A_j$ and $|\gamma|\lesssim 2^{-j(d-1)/2}$ to
ensure the functions are positive. The constant is adjusted so that
for one of the $f_m$ that we denote $f_0$, $\forall x\in H^+,\
\left|\cH(f_0^-)(x)\right|\le c_b$
with $c_b\in(0,\frac12)$. The function  $f_m$ also integrate to 1. We denote by
$\mathcal{M}$ the cardinality of $A_j$ ($\mathcal{M}\simeq2^{j(d-1)}$), $P_{m}$ the distributions of an i.i.d. sample of $(Y,X)$ of size $n$ when $\fbeta = f_m$ and for a given $\fX$, and $\Lambda(P_m,P_0)$ the likelihood ratio. Recall that
$K(P_m,P_0)=\Esp_{P_m}\left[\Lambda(P_m,P_0)\right]$. 
We make use of the following Lemma from \cite{KT}. 
\begin{lemma}\label{lem6}
If for $\pi_0>0$ and  $\mathcal{M}\in\xN^*$ the following three condition hold
\begin{enumerate}[\textup{(}i\textup{)}]
\item\label{lem6i}
$f_m\in B_{r,q}^s(M)\cap\mathcal{D}$ for $m=1,\hdots,\mathcal{M}$,
\item\label{lem6ii} $\forall m\ne l,\ \|f_m-f_l\|_p\ge 2h>0$,
\item\label{lem6iii} $\forall m=1,\hdots,\mathcal{M}$, $\Lambda(P_0,P_m)=\exp(z_n^m-v_n^m)$, where $z_n^m$
are random variables and $v_n^m$ constants such that
$\mathbb{P}(z_n^m>0)\ge\pi_0$ and
$\exp\left(\sup_{m=1,\hdots,\mathcal{M}}v_n^m\right)\le
\mathcal{M}$,
\end{enumerate}
then
\[
\forall z\ge1,\ \inf_{\hatfbeta}\sup_{\fbeta\in B^s_{r,q}(M)\cap\mathcal{D}}
\Esp\left\|\hatfbeta-\fbeta\right\|_p^z \ge \frac{h^{-z}\pi_0}{2}.
\]
\end{lemma}
\noindent Item \eqref{lem6i} is satisfied when $|\gamma|\le
M2^{-j(s-(d-1)(1/r-1/2)}$. This is more restrictive than the
condition to ensure positivity because we assume that $s\ge(d-1)/r$.
Thus, now we take $\gamma= 2cM 2^{-j(s-(d-1)(1/r-1/2)}$ for a well-chosen constant $c$.\\ 
The constant $h$ in
\eqref{lem6ii} is obtained as follows, if $m\neq m'$,
\begin{align*}
\|f_m-f_{m'}\|_p&=\gamma\|\psi_{j,\xi_m}-\psi_{j,\xi_{m'}}\|_p\\
&\ge \gamma c_{p,A} 2^{j(d-1)(1/2-1/p)}\\
&\ge 2c M 2^{-j(s-(d-1)(1/r-1/p))}.
\end{align*}
Let us now consider item \eqref{lem6iii}, we obtain
\begin{align*}
P_m\left(\log\left(\Lambda(P_0,P_m)\right)\ge-j(d-1)\log 2\right)
&\ge1-P_m\left(\left|\log\left(\Lambda(P_0,P_m)\right)\right|\ge j(d-1)\log 2\right)\\
&\ge1-\frac{\Esp_{P_m}\left[\left|\log\left(\Lambda(P_0,P_m)\right)\right|\right]}{j(d-1)\log 2}.
\end{align*}
Thus, condition \eqref{lem6iii} is satisfied when
\[
\Esp_{P_m}\left[\left|\log\left(\Lambda(P_0,P_m)\right)\right|\right]\le
\alpha j(d-1)\log 2,
\]
for $\alpha\in(0,1)$.  The same computations as in the beginning of
Section \ref{s41} yield that we need to impose
$n2^{-2j\nu(d)}\gamma^2\lesssim j$, thus
\[
A_Xn2^{-2j(s+\nu(d)-(d-1)(1/r-1/2))}\lesssim j.
\]
The desired rate is obtained by taking
\[
2^j\simeq\left(\frac{nA_X}{\log
    \left(nA_X\right)}\right)^{\frac{1}{2(s+\nu(d)-(d-1)(1/r-1/2))}}.
\]

\subsection{Comparison between Besov ellipsoids of a function and its odd part}
\begin{lemma}\label{compnorms}
For $0<s,q\le\infty$ and $1\le r\le\infty$, there exists a constant $c_{\mathrm{eq}}$ that can depend on $d$ such that, for every $f\in B^s_{r,q}$, 
 $\|f^-\|_{B^s_{r,q}}\le c_{\mathrm{eq}}\|f\|_{B^s_{r,q}}$.
\end{lemma}
\noindent {\bf Proof.} In Definition \ref{def2} every $f\in B^s_{r,q}(\xSd)$ has same norm as $x\to f(-x)$, thus by the triangle inequality $\|f^-\|_{B^s_{r,q}}^A\le \|f\|_{B^s_{r,q}}^A$. We conclude by equivalence of the norms.\hfill $\square$

\subsection{A general inequality}\label{s43}
We make use of the constants $c_{1,z}$ and $c_{2,z}$ such that
\begin{align}
  \int_{\bR^+} z \tau^{z-1} e^{-\beta\tau} d\tau &\leq c_{1,z}
  \beta^{-z}\label{UBexp}\\
  \int_{\bR^+} z \tau^{z-1} e^{-\alpha\tau^2} d\tau &\leq c_{2,z}
  \alpha^{-z/2}\label{UBnorm}.
\end{align}
\begin{lemma}\label{toracle}
For every $\tau,\ \gamma,\ z>1$ and
\[
T_{j,\xi,\gamma}^{s,++}\ge 3\sqrt{2\gamma}t_n\widehat{\sigma}_{j,\xi}+26M_{j,\xi}\frac{\gamma \log n}{n-1}\myeq T_{j,\xi,\gamma}^{s,+},
\] the two following inequalities hold:\\
when $p=\infty$,
\begin{align*}
&\frac{1}{2^{z-1}}
\mathbb{E}\left[\left\|\hatfbeta^{a,\rho} -\fbeta\right\|_{\infty}^z\right]\\
&\leq
\left\|\fbetam^{a,J} -\fbetam\right\|_{\infty}^z+ (J+1)^{z-1}\ConstNeedEquiv^z_{\infty}\Big\{\\
&\qquad a_{n,\infty,z,J}
\sum_ {j=0}^J
 2^{j(d-1)z/2}\left(\sup_{\xi \in \Quad_j}
\left|\beta^{a}_{j,\xi}\right|^z
\ind_{\left|\beta^{a}_{j,\xi}\right| \leq T^{s,++}_{j,\xi,\gamma} } + \Esp\left[ \sup_{\xi \in \Quad_j}\left|\hatbeta^a_{j,\xi}-\beta^{a}_{j,\xi}\right|^z\ind_{\left|\beta^{a}_{j,\xi}\right|> T^{s,++}_{j,\xi,\gamma}}\right]
\right)\\
& \qquad+
\frac{4C_{\Xi}}{n^\gamma}
  \sum_{j=0}^J
2^{j(d-1)(z/2+1)}
\sup_{\xi\in\Quad_j}
  \left|
    \beta^{a}_{j,\xi}\right|^z\\
&\qquad + \left(\frac{C_{\Xi}4}{n^{\gamma}}\right)^{1-1/\tau}\left(\frac{1}{\sqrt{n}}B_X^{1/2}2^{Jz(\nu(d)+(d-1)/2)}\right)^{z}
2^{J(d-1)(1-1/\tau)}b_{n,\infty,z,J,\tau}\Big\},
\end{align*}
where
\begin{align*}
a_{n,\infty,z,J}&=1+
\left(\frac{2}{\sqrt{\gamma \log n}} \right)^z  \left(2+
\left(\log\left(C_{\Xi}2^{J(d-1)}c_{2,z}\right)\right)^{z/2}\right)\\
&\quad+\left(\frac{4}{\gamma \log n}\right)^z \left(2+
\left(\log\left(C_{\Xi}2^{J(d-1)}c_{1,z}\right)\right)^{z}\right)\\
b_{n,\infty,z,J,\tau}&=
\frac{\left(2\sqrt{2}C_2B(d,2)\right)^z\left(2^{1/\tau}+
\left(\log\left(C_{\Xi}2^{J(d-1)}c_{2,z\tau}\right)\right)^{z/2}\right)}
{1-2^{-(z\nu(d)+(d-1)(z/2+1-1/\tau))}}\\
&\quad+\frac{\left(8C_{\infty}B(d,\infty)/3\right)^z\left(2^{1/\tau}+
\left(\log\left(C_{\Xi}2^{J(d-1)}c_{1,z}\right)\right)^{z}\right)}
{1-2^{-(z\nu(d)+(d-1)(z+1-1/\tau))}}
\left(\frac{2^{J(d-1)}}{n}B_X\right)^{z/2};
\end{align*}
while, when $1\le p\le \infty$,
\begin{align*}
&\frac{1}{2^{z-1}}
\mathbb{E}\left[\left\|\hatfbeta^{a,\rho} -\fbeta\right\|_p^z\right]\\
&\leq
\left\|\fbetam^{a,J} -\fbetam\right\|_p^z+ (J+1)^{z-1}\ConstNeedEquiv^z_pC_{\Xi}^{z/(p\wedge z)-1}\Big\{\\
&\qquad a_{n,p,z,J}\sum_{j=0}^J
 2^{j(d-1)(z/2-z/(p\vee z))}\sum_{\xi\in\Quad_j}\left(
\left|\beta^{a}_{j,\xi}\right|^z
\ind_{\left|\beta^{a}_{j,\xi}\right| \leq T^{s,++}_{j,\xi,\gamma} } + \Esp\left[ \left|\hatbeta^a_{j,\xi}-\beta^{a}_{j,\xi}\right|^z\right]
\ind_{\left|\beta^{a}_{j,\xi}\right|> T^{s,++}_{j,\xi,\gamma}}\right)\\
&\qquad+ \frac{4}{n^\gamma}\sum_ {j=0}^J
 2^{j(d-1)z(1/2-1/(p\vee z))} \sum_{\xi\in\Quad_j}
  \left|
    \beta^{
    a}_{j,\xi}\right|^z\\
&\qquad + \frac{2^{2-1/\tau}}{n^{\gamma(1-\frac{1}{\tau})}}C_{\Xi}\left(\frac{1}{\sqrt{n}}B_X^{1/2}2^{J(\nu(d)+(d-1)/2)}\right)^{z}
2^{J(d-1)(1-z/(p\vee z))}b_{n,p,z,J,\tau}\Big\},
\end{align*}
where
\begin{align*}
a_{n,p,z,J}&=1+2 \left( \left(
\frac{\sqrt{2} c_{2,z}^{1/z}}{\sqrt{\gamma\log n}}
\right)^z
+
\left(
\frac{2 c_{1,z}^{1/z}}{\gamma \log n}
\right)^z
\right)\\
b_{n,p,z,J,\tau}&=\frac{\left( 2 c_{2,z\tau}^{1/(z\tau)}
\ConstNeedNorm_{2}B(d,2) \right)^z}{1-2^{-(z\nu(d)+(d-1)(z/2+1-z/(p\vee z)))}}
+ \frac{\left(\frac{4}{3} c_{1,z\tau}^{1/(z\tau)}
\ConstNeedNorm_{\infty}B(d,\infty)\right)^z}
{1-2^{-(z\nu(d)+(d-1)(z+1-z/(p\vee z)))}}
\left(\frac{2^{J(d-1)}}{n}B_X\right)^{z/2}.
\end{align*}
\end{lemma}
The inequalities of Lemma \ref{toracle} are similar to oracle inequalities, for a well-chosen $J$ depending on $n$ (see Theorem \ref{t2}),
where the oracle
estimates $\beta^{a}_{j,\xi}$ if and only if the error made by
estimating this coefficient is smaller than the one made by
discarding it.  This \emph{oracle} strategy would lead to a quantity
of the form
\begin{align*}
\left|
    \beta^{a}_{j,\xi}\right|^z
\ind_{\left|\beta^{a}_{j,\xi}\right|\leq
    \left(\Esp\left[ \left|\hatbeta^a_{j,\xi}-\beta^{a}_{j,\xi}\right|^z\right]\right)^{1/z}}
+
\Esp\left[ \left|\hatbeta^a_{j,\xi}-\beta^{a}_{j,\xi}\right|^z\right]
\ind_{\left|\beta^{a}_{j,\xi}\right|>
    \left(\Esp\left[ \left|\hatbeta^a_{j,\xi}-\beta^{a}_{j,\xi}\right|^z\right]\right)^{1/z}}.
\end{align*}
Proving such an oracle inequality would require to lower bound
$\left(\Esp\left[
    \left|\hatbeta^a_{j,\xi}-\beta^{a}_{j,\xi}\right|^z\right]\right)^{1/z}$.
In the inequalities of Lemma \ref{toracle} the ideal quantity
$\left(\Esp\left[ \left|\hatbeta^a_{j,\xi}-\beta^{a}_{j,\xi}\right|^z\right]\right)^{1/z}$
is replaced by $T_{j,\xi,\gamma}^{s,++}$, called {\em quasi-oracle}.
The remaining terms can be made as
small as we want by taking $\gamma$ large enough.
The last term corresponds to the approximation error.
Upper bounds of these types, uniform on Besov
ellipsoids, yield an approximation error which can
be expressed in terms of the regularity of the Besov class
and is uniformly small for $J$ large enough and allows to treat the bias/variance trade-off
in the quasi-oracle term uniformly over the ellipsoid.

\subsection{Proof of Lemma \ref{toracle}}\label{s7}

\subsubsection{Preliminaries}\label{s71}
Recall from the proof of Theorem 4.1 in \cite{GK} that for every $1\le p\le \infty$
\[
\left\|\hatfbeta^{a,\rho} -\fbeta\right\|_p\le2\left\|\hatfbetam^{a,\rho} 
  - \fbetam\right\|_p,
\]
and that, for
$1\le z< \infty$, we have
\begin{equation}\label{UB1}
\left\|\hatfbetam^{a,\rho} 
-\fbetam\right\|_p^z\le2^{z-1}\left(\left\|\hatfbetam^{a,\rho} -\fbetam^{a,J}\right\|_p^z
+\left\|\fbetam^{a,J} -\fbetam\right\|_p^z\right).
\end{equation}
The first term corresponds to the error in the high dimensional
space while the second term corresponds to the approximation error.  Let us start by studying the first term.\\
Lemma \ref{lem1} \eqref{lem1i} yields
\begin{align*}
  \left\| \hatfbetam^{a,\rho} -\fbetam^{a,J} \right\|_p^z
&\leq (J+1)^{z-1} \sum_{j=0}^J
\left\| \sum_{\xi\in\Quad_j}
\left( \rho_{T_{j,\xi,\gamma}}\left(\hatbeta^a_{j,\xi}\right) - \beta^{a}_{j,\xi} \right)
\psi_{j,\xi} \right\|_p^z\\
&\leq (J+1)^{z-1} \sum_{j=0}^J
\ConstNeedEquiv^z_p 2^{j(d-1)z(1/2-1/p)}
\left\|
\rho_{T_{j,\xi,\gamma}}\left(\hatbeta^a_{j,\xi}\right) - \beta^{a}_{j,\xi}
\right\|_p^z.
\end{align*}
Thus, for $p=\infty$, we have
\begin{align*}
  \left\| \hatfbetam^{a,\rho} -\fbetam^{a,J} \right\|_p^z
& \leq (J+1)^{z-1} \sum_{j=0}^J
\ConstNeedEquiv^z_{\infty} 2^{j(d-1)z/2} \sup_{\xi\in\Quad_j} \left|
\rho_{T_{j,\xi,\gamma}}\left(\hatbeta^a_{j,\xi}\right) - \beta^{a}_{j,\xi}
\right|^z,
&\end{align*}
while, for $p<\infty$, we have
\begin{align*}
  \left\| \hatfbetam^{a,\rho} -\fbetam^{a,J} \right\|_p^z
&\leq (J+1)^{z-1}\ConstNeedEquiv^z_p C_{\Xi}^{z/(p\wedge z)-1}\sum_ {j=0}^J
 2^{j(d-1)z(1/2-1/(p\vee z))} \sum_{\xi\in\Quad_j} \left|\rho_{T_{j,\xi,\gamma}}\left(\hatbeta^a_{j,\xi}\right) - \beta^{a}_{j,\xi}\right|^z.
\end{align*}
The last inequality is obtained by using that, when $p\ge z$, we have $$
\left(\sum_{\xi\in\Quad_j}
\left|b_{\xi}\right|^p\right)^{z/p}\le\sum_{\xi\in\Quad_j}
\left|b_{\xi}\right|^z,$$ and by the H\"older inequality, when $p\le z$, we have
\[
\left(\sum_{\xi\in\Quad_j}
\left|b_{\xi}\right|^p\right)^{z/p}\le
C_{\Xi}^{z/p-1}\sum_{\xi\in\Quad_j} \left|b_{\xi}\right|^z.
\]

\subsubsection{Coefficientwise analysis}
For the simplicity of the notations we sometimes drop the dependence on $\gamma$ in the sets of indices.\\
We first consider the term
\begin{align*}
\delta_{j,\xi,z} &\myeq \left|\rho_{T_{j,\xi,\gamma}}\left(\hatbeta^a_{j,\xi}\right) - \beta^{a}_{j,\xi}\right|^z.
\end{align*}
By construction we have
\begin{align*}
 \delta_{j,\xi,z}
& =    \left|
    \beta^{a}_{j,\xi}\right|^z
  \ind_{\left|\hatbeta^a_{j,\xi}\right|\leq T_{j,\xi,\gamma}}
+   \left|\hatbeta^a_{j,\xi} -
    \beta^{a}_{j,\xi}\right|^z
  \ind_{\left|\hatbeta^a_{j,\xi}\right|> T_{j,\xi,\gamma}}\\
& = \max \left( \left|
    \beta^{a}_{j,\xi}\right|^z
  \ind_{\left|\hatbeta^a_{j,\xi}\right|\leq T_{j,\xi,\gamma}}
,  \left|\hatbeta^a_{j,\xi} -
    \beta^{a}_{j,\xi}\right|^z
  \ind_{\left|\hatbeta^a_{j,\xi}\right|> T_{j,\xi,\gamma}} \right).
\end{align*}
We introduce two ``phantom'' random thresholds
$T^{b}_{j,\xi,\gamma}= T_{j,\xi,\gamma} - \Delta_{j,\xi,\gamma}$ and
$T^{s}_{j,\xi,\gamma}=T_{j,\xi,\gamma} + \Delta_{j,\xi,\gamma}$ for
some $\Delta_{j,\xi,\gamma}$ to be defined later.  They are used to
define \emph{big} and \emph{small} original needlet coefficients. We
also use $T^{b,-}_{j,\xi,\gamma}$ for a deterministic 
lower bound on $T^{b}_{j,\xi,\gamma}$, $T^{s,+}_{j,\xi,\gamma}$ and
$\Delta^{+}_{j,\xi,\gamma}$ for deterministic upper bounds on $T^{s}_{j,\xi,\gamma}$ and 
$\Delta_{j,\xi,\gamma}$. These bounds will hold with high probability.  We obtain almost surely
\begin{align*}
 \delta_{j,\xi,z}
& =   \max \bigg( \left|
    \beta^{a}_{j,\xi}\right|^z
\max\left(
  \ind_{\left|\hatbeta^a_{j,\xi}\right|\leq T_{j,\xi,\gamma}}
\ind_{\left|\beta^{a}_{j,\xi}\right|\leq T^{s}_{j,\xi,\gamma}}
,
  \ind_{\left|\hatbeta^a_{j,\xi}\right|\leq T_{j,\xi,\gamma}}
\ind_{\left|\beta^{a}_{j,\xi}\right|> T^{s}_{j,\xi,\gamma}}
\right),\\
&\qquad\qquad
\left|\hatbeta^a_{j,\xi} -
    \beta^{a}_{j,\xi}\right|^z
\max\left(
  \ind_{\left|\hatbeta^a_{j,\xi}\right|> T_{j,\xi,\gamma}} \ind_{\left|\beta^{a}_{j,\xi}\right|\leq T^{b}_{j,\xi,\gamma}}
,
  \ind_{\left|\hatbeta^a_{j,\xi}\right|> T_{j,\xi,\gamma}} \ind_{\left|\beta^{a}_{j,\xi}\right|> T^{b}_{j,\xi,\gamma}}
\right)
\bigg)\\
& \leq
\max \bigg( \left|
    \beta^{a}_{j,\xi}\right|^z
\max\left(
\ind_{\left|\beta^{a}_{j,\xi}\right|\leq T^{s}_{j,\xi,\gamma}}
,
  \ind_{\left|\hatbeta^a_{j,\xi}-\beta^{a}_{j,\xi}\right|>\Delta_{j,\xi,\gamma}}
\right),\\
&\qquad\qquad\left|\hatbeta^a_{j,\xi} -
    \beta^{a}_{j,\xi}\right|^z
\max\left(
  \ind_{\left|\hatbeta^a_{j,\xi}-\beta^{a}_{j,\xi}\right|>\Delta_{j,\xi,\gamma}}
,
 \ind_{\left|\beta^{a}_{j,\xi}\right|> T^{b}_{j,\xi,\gamma}}
\right)
\bigg)\\
& \leq
\max \bigg( \left|
    \beta^{a}_{j,\xi}\right|^z
\max\left(
\ind_{\left|\beta^{a}_{j,\xi}\right|\leq
    T^{s,+}_{j,\xi,\gamma}}
,
\ind_{T^{s,+}_{j,\xi,\gamma} < T^{s}_{j,\xi,\gamma}}
,
  \ind_{\left|\hatbeta^a_{j,\xi}-\beta^{a}_{j,\xi}\right|>\Delta_{j,\xi,\gamma}}
\right),\\
&\qquad\qquad
\left|\hatbeta^a_{j,\xi} -
    \beta^{a}_{j,\xi}\right|^z
\max\left(
  \ind_{\left|\hatbeta^a_{j,\xi}-\beta^{a}_{j,\xi}\right|>\Delta_{j,\xi,\gamma}}
,
\ind_{\left|\beta^{a}_{j,\xi}\right|> T^{b,-}_{j,\xi,\gamma}},
\ind_{ T^{b,-}_{j,\xi,\gamma}> T^{b}_{j,\xi,\gamma}}
\right)
\bigg).
\end{align*}
Sorting the terms according to the number of random terms we obtain
\begin{align*}
\delta_{j,\xi,z}& \leq
\max \bigg(
\left|
    \beta^{a}_{j,\xi}\right|^z
\ind_{\left|\beta^{a}_{j,\xi}\right|\leq
    T^{s,+}_{j,\xi,\gamma}}
,
\left|
    \beta^{a}_{j,\xi}\right|^z
\max\left(
\ind_{T^{s,+}_{j,\xi,\gamma} < T^{s}_{j,\xi,\gamma}}
,
  \ind_{\left|\hatbeta^a_{j,\xi}-\beta^{a}_{j,\xi}\right|>\Delta_{j,\xi,\gamma}}
\right),\\
&\qquad\qquad
  \left|\hatbeta^a_{j,\xi} -
    \beta^{a}_{j,\xi}\right|^z
\ind_{\left|\beta^{a}_{j,\xi}\right|> T^{b,-}_{j,\xi,\gamma}}
,
\left|\hatbeta^a_{j,\xi} -
    \beta^{a}_{j,\xi}\right|^z
\max\left(
  \ind_{\left|\hatbeta^a_{j,\xi}-\beta^{a}_{j,\xi}\right|>\Delta_{j,\xi,\gamma}}
,
\ind_{ T^{b,-}_{j,\xi,\gamma}> T^{b}_{j,\xi,\gamma}}
\right)
\bigg).
\end{align*}

\subsubsection{Scalewise analysis}

Defining 
\begin{align*}
  M_{j,z} &\myeq \sup_{\xi\in\Quad_j} \left|
\rho_{T_{j,\xi,\gamma}}\left(\hatbeta^a_{j,\xi}\right) - \beta^{a}_{j,\xi}
\right|^z = \sup_{\xi\in\Quad_j} \delta_{j,\xi,z}\\
S_{j,z} & \myeq\sum_{\xi\in\Quad_j}
  \left|\rho_{T_{j,\xi,\gamma}}\left(\hatbeta^a_{j,\xi}\right) -
    \beta^{a}_{j,\xi}\right|^z
= \sum_{\xi\in\Quad_j} \delta_{j,\xi,z},
\end{align*}
we obtain
\begin{align*}
M_{j,z} &\leq \max\bigg( \sup_{\xi\in\Quad_j}
  \left|
    \beta^{a}_{j,\xi}\right|^z
\ind_{\left|\beta^{a}_{j,\xi}\right|\leq
    T^{s,+}_{j,\xi,\gamma}},
\sup_{\xi\in\Quad_j}
  \left|
    \beta^{a}_{j,\xi}\right|^z
\max \left(
\ind_{T^{s,+}_{j,\xi,\gamma} < T^{s}_{j,\xi,\gamma}},
  \ind_{\left|\hatbeta^a_{j,\xi}-\beta^{a}_{j,\xi}\right|>\Delta_{j,\xi,\gamma}}
\right),
\\
& \quad
\sup_{\xi\in\Quad_j} \left|\hatbeta^a_{j,\xi} -
    \beta^{a}_{j,\xi}\right|^z
 \ind_{\left|\beta^{a}_{j,\xi}\right|> T^{b,-}_{j,\xi,\gamma}}
,
\sup_{\xi\in\Quad_j} \left|\hatbeta^a_{j,\xi} -
    \beta^{a}_{j,\xi}\right|^z
\max\left(
\ind_{ T^{b,-}_{j,\xi,\gamma}> T^{b}_{j,\xi,\gamma}}
,   \ind_{\left|\hatbeta^a_{j,\xi}-\beta^{a}_{j,\xi}\right|>
    \Delta_{j,\xi,\gamma}}
\right)\bigg) \\
&\myeq \max(M_{j,z}^{S0},M_{j,z}^{S1},M_{j,z}^{B1},M_{j,z}^{B2})\\
&\leq M_{j,z}^{S0} + M_{j,z}^{S1} + M_{j,z}^{B1} + M_{j,z}^{B2};
\\
S_{j,z}& \leq \sum_{\xi\in\Quad_j}
  \left|
    \beta^{a}_{j,\xi}\right|^z
\ind_{\left|\beta^{a}_{j,\xi}\right|\leq
    T^{s,+}_{j,\xi,\gamma}}
+
\sum_{\xi\in\Quad_j}
  \left|
    \beta^{a}_{j,\xi}\right|^z
\max\left(
\ind_{T^{s,+}_{j,\xi,\gamma} < T^{s}_{j,\xi,\gamma}}
,
  \ind_{\left|\hatbeta^a_{j,\xi}-\beta^{a}_{j,\xi}\right|>\Delta_{j,\xi,\gamma}}
\right)
\\
& \quad +
\sum_{\xi\in\Quad_j} \left|\hatbeta^a_{j,\xi} -
    \beta^{a}_{j,\xi}\right|^z
\ind_{\left|\beta^{a}_{j,\xi}\right|> T^{b,-}_{j,\xi,\gamma}}
+
\sum_{\xi\in\Quad_j} \left|\hatbeta^a_{j,\xi} -
    \beta^{a}_{j,\xi}\right|^z
\max\left(
\ind_{ T^{b,-}_{j,\xi,\gamma}> T^{b}_{j,\xi,\gamma}}
,   \ind_{\left|\hatbeta^a_{j,\xi}-\beta^{a}_{j,\xi}\right|>
    \Delta_{j,\xi,\gamma}}
\right) \\
&\myeq S_{j,z}^{S0}+S_{j,z}^{S1}+S_{j,z}^{B1}+S_{j,z}^{B2}.
\end{align*}
We bound the expectations of the random terms as follows
\begin{align*}
\Esp\left[ M_{j,z}^{S1} \right] & \leq
\sup_{\xi\in\Quad_j}
  \left|
    \beta^{a}_{j,\xi}\right|^z
\Esp \left[
\sup_{\xi\in\Quad_j}
\max \left(
\ind_{T^{s,+}_{j,\xi,\gamma} < T^{s}_{j,\xi,\gamma}},
  \ind_{\left|\hatbeta^a_{j,\xi}-\beta^{a}_{j,\xi}\right|>\Delta_{j,\xi,\gamma}}
\right)
\right]\\
& \leq \sup_{\xi\in\Quad_j}
  \left|
    \beta^{a}_{j,\xi}\right|^z
\left(
\Prob\left(
\bigcup_{\xi \in \Quad_j}\left\{T^{s,+}_{j,\xi,\gamma} < T^{s}_{j,\xi,\gamma}\right\}\right)+
  \Prob\left(
\bigcup_{\xi \in \Quad_j}\left\{\left|\hatbeta^a_{j,\xi}-\beta^{a}_{j,\xi}\right|>\Delta_{j,\xi,\gamma}\right\}\right)
\right);
\\
\Esp\left[ M_{j,z}^{B1} \right] & \leq
\Esp \left[ \sup_{\xi\in\Quad_j} \left|\hatbeta^a_{j,\xi} -
    \beta^{a}_{j,\xi}\right|^z
 \ind_{\left|\beta^{a}_{j,\xi}\right|> T^{b,-}_{j,\xi,\gamma}}
\right];\\
\Esp\left[ M_{j,z}^{B2} \right] & \leq
\Esp \left[ \sup_{\xi\in\Quad_j} \left|\hatbeta^a_{j,\xi} -
    \beta^{a}_{j,\xi}\right|^{z\tau}
\right]^{1/\tau}
\Esp \left[
\sup_{\xi\in\Quad_j}
\max\left(
\ind_{ T^{b,-}_{j,\xi,\gamma}> T^{b}_{j,\xi,\gamma}}
,   \ind_{\left|\hatbeta^a_{j,\xi}-\beta^{a}_{j,\xi}\right|>
    \Delta_{j,\xi,\gamma}}
\right)\right]^{1-1/\tau}\\
& \leq
\Esp \left[ \sup_{\xi\in\Quad_j} \left|\hatbeta^a_{j,\xi} -
    \beta^{a}_{j,\xi}\right|^{z\tau}
\right]^{1/\tau}
\left(
\Prob\left(
\bigcup_{\xi \in \Quad_j}\left\{T^{b,-}_{j,\xi,\gamma}> T^{b}_{j,\xi,\gamma}
\right\}\right)
+
\Prob\left(
\bigcup_{\xi \in \Quad_j}\left\{\left|\hatbeta^a_{j,\xi}-\beta^{a}_{j,\xi}\right|>
    \Delta_{j,\xi,\gamma}
\right\}\right)
\right)^{1-1/\tau};\\
  \Esp\left[ S_{j,z}^{S1} \right] &= \sum_{\xi\in\Quad_j}
  \left|
    \beta^{a}_{j,\xi}\right|^z
 \Esp \left[\max\left(
\ind_{T^{s,+}_{j,\xi,\gamma} < T^{s}_{j,\xi,\gamma}}
,\ind_{\left|\hatbeta^a_{j,\xi}-\beta^{a}_{j,\xi}\right|>\Delta_{j,\xi,\gamma}}
\right) \right]\\
& \le \sum_{\xi\in\Quad_j}
  \left|
    \beta^{a}_{j,\xi}\right|^z \left(
\Prob\left\{ T^{s,+}_{j,\xi,\gamma} < T^{s}_{j,\xi,\gamma} \right\}
+
\Prob\left\{ \left|\hatbeta^a_{j,\xi}-\beta^{a}_{j,\xi}\right|>\Delta_{j,\xi,\gamma}\right\}
\right);\\
  \Esp\left[ S_{j,z}^{B1} \right] &= \sum_{\xi\in\Quad_j} \Esp\left[\left|\hatbeta^a_{j,\xi} -
    \beta^{a}_{j,\xi}\right|^z \right]
\ind_{\left|\beta^{a}_{j,\xi}\right|> T^{b,-}_{j,\xi,\gamma}};\\
  \Esp\left[ S_{j,z}^{B2} \right] &= \sum_{\xi\in\Quad_j} \Esp \left[ \left|\hatbeta^a_{j,\xi} -
   \beta^{a}_{j,\xi}\right|^z
\max\left(
\ind_{ T^{b,-}_{j,\xi,\gamma}> T^{b}_{j,\xi,\gamma}}
,   \ind_{\left|\hatbeta^a_{j,\xi}-\beta^{a}_{j,\xi}\right|>
    \Delta_{j,\xi,\gamma}}
\right) \right]\\
& \leq \sum_{\xi\in\Quad_j} \left( \Esp \left[ \left|\hatbeta^a_{j,\xi} -
    \beta^{a}_{j,\xi}\right|^{z\tau} \right]\right)^{1/\tau}
\left(
\Prob\left\{ T^{b,-}_{j,\xi,\gamma}> T^{b}_{j,\xi,\gamma}\right\}
+
\Prob\left\{\left|\hatbeta^a_{j,\xi}-\beta^{a}_{j,\xi}\right|>
    \Delta_{j,\xi,\gamma}\right\}
\right)^{1-1/\tau}.
\end{align*}
The constant $\tau>1$ in the H\"older inequality will be specified later.\\

\subsubsection{Bernstein inequality and the term $ \left|\hatbeta^a_{j,\xi} -
      \beta^{a}_{j,\xi}\right|^z$}\label{s721}
Let denote variance of $G_{j,\xi}(X,Y)$
\begin{align*}
\sigma_{j,\xi}^2 \myeq\Esp\left[ \left( G_{j,\xi}(X,Y) - \beta^{a}_{j,\xi}\right)^2 \right].
\end{align*}
\begin{lemma}\label{lem7}
We have
\begin{equation*}
  \Esp \left[
\left|\hatbeta^a_{j,\xi} -\beta^{a}_{j,\xi} \right|
^z
\right]\le
 2 \left(
c_{2,z}
\left(
\frac{2}{\sqrt{n}} \sigma_{j,\xi}
\right)^z
+
c_{1,z}
\left(
\frac{4}{3 n}
M_{j,\xi}
\right)^z
\right).
\end{equation*}
\end{lemma}
\noindent {\bf Proof.}
The Bernstein inequality yields
\begin{align*}
  \Prob\left\{ \left|\hatbeta^a_{j,\xi} -\beta^{a}_{j,\xi} \right| \geq u
  \right\}
&\leq 2 e^{-\frac{nu^2}{2\left(\left(\sigma_{j,\xi}\right)^2 + M_{j,\xi} u /3\right)}}\\
&\le  2 \left(e^{-\frac{nu^2}{4\left(\sigma_{j,\xi}\right)^2}}+e^{-\frac{3nu}{4M_{j,\xi}}}\right).
\end{align*}
Using now
$\Esp\left[ |X|^z \right]  = \int_{\bR^+} z u^{z-1} \Prob\{ |X|>u\}du$, we obtain
\begin{align*}
  \Esp \left[
\left|\hatbeta^a_{j,\xi} -\beta^{a}_{j,\xi} \right|
^z
\right]
& \leq \int_{\bR^+} z u^{z-1} \Prob\left\{ \left|\hatbeta^a_{j,\xi} -\beta^{a}_{j,\xi} \right| \geq u
  \right\} du\\
& \leq \int_{\bR^+} z u^{z-1} 2
  \left(e^{-\frac{nu^2}{4\left(\sigma_{j,\xi}\right)^2}}+e^{-\frac{3nu}{4M_{j,\xi}}}\right)
  du,
\end{align*}
hence the inequality from the lemma follows from \eqref{UBexp} and \eqref{UBnorm}.\hfill $\square$\vspace{0.3cm}

Lemma \ref{lem7} is used to obtain a uniform upper bound of the power
of the ratio between $\left|\hatbeta^a_{j,\xi} -\beta^{a}_{j,\xi}
\right|$ and a threshold $c_{\sigma} \sqrt{\log(n)/n} \sigma_{j,\xi} + c_{M}
\log(n)/(n-1) M_{j,\xi}$:
\begin{align}
 \Esp \left[
\left(\frac{
\left|\hatbeta^a_{j,\xi} -\beta^{a}_{j,\xi} \right|
}{
c_\sigma \sqrt{\log(n)/n}
\sigma_{j,\xi} + c_M \log(n)/(n-1) M_{j,\xi}
}\right)^z
\right]
& \leq
2 \left(
c_{2,z}
\left(
2 \frac{1}{c_\sigma \sqrt{\log n}+ c_M \frac{\sqrt{n}\log n}{n-1} \frac{M_{j,\xi}}{\sigma_{j,\xi}}}
\right)^z\right.\nonumber\\
&\quad
\left.+
c_{1,z}
\left(
\frac{4}{3}
\frac{1}{c_\sigma \sqrt{n} \sqrt{\log n} \frac{\sigma_{j,\xi}}{M_{j,\xi}} +
  c_M \log(n) \frac{n}{n-1}  }
\right)^z
\right)\nonumber\\
& \leq
2 \left(
c_{2,z}
\left(
2 \frac{1}{c_\sigma \sqrt{\log n}}
\right)^z
+
c_{1,z}
\left(
\frac{4}{3}
\frac{1}{c_M\log n  }
\right)^z
\right).\label{elem7ii}
\end{align}
The following similar lemma is useful to handle the case $p=\infty$.
\begin{lemma}\label{lem8}
For any $\Xi'_j \subset \Xi_j$, we have
\begin{align}
 \Esp\left[ \sup_{\xi \in \Xi'_j}  \left(\frac{\left|\hatbeta^a_{j,\xi} -\beta^{a}_{j,\xi}
   \right|}{u_{j,\xi}}\right)^z
\right]& \leq   \left(\frac{2\sqrt{2}}{\sqrt{n}} \sup_{\xi \in \Xi'_j} \frac{\sigma_{j,\xi}}{c_{j,\xi}}\right)^z  \left(2+
\left(\log\left(c_{2,z} \left|\Xi'_j\right|\right)\right)^{z/2} \right)\nonumber\\
&\quad
+ \left(\frac{8}{3n} \sup_{\xi \in \Xi'_j}
  \frac{M_{j,\xi}}{c_{j,\xi}} \right)^z \left(2+\left(\log\left( c_{1,z}\left|\Xi'_j\right|\right)\right)^z \right).\label{elem8}
\end{align}
\end{lemma}
\noindent {\bf Proof.} A uniform union bound yields
\begin{align*}
&\Prob\left\{ \sup_{\xi \in \Xi'_j}  \frac{\left|\hatbeta^a_{j,\xi} -\beta^{a}_{j,\xi}
   \right|}{u_{j,\xi}} \geq \tau
\right\}\\
& \leq \min\left(1,\left|\Xi'_j\right| 2 \left(
e^{-\frac{1}{4} n  \left(\inf_{\xi \in \Xi'_{j}} \frac{u_{j,\xi}}{\sigma_{j,\xi}}  \right)^2\tau^2 }
+
e^{-\frac{3}{4}n \left(\inf_{\xi \in \Xi'_{j}}
  \frac{u_{j,\xi}}{M_{j,\xi}}  \right) \tau } \right)\right)
\\
& \leq \min\left(1,\left|\Xi'_j\right| 2 
e^{-\frac{1}{4} n  \left(\inf_{\xi \in \Xi'_{j}} \frac{u_{j,\xi}}{\sigma_{j,\xi}}  \right)^2\tau^2 }
\right)
+ \min\left(1,\left|\Xi'_j\right| 2 
e^{-\frac{3}{4}n \left(\inf_{\xi \in \Xi'_{j}}
  \frac{u_{j,\xi}}{M_{j,\xi}}  \right) \tau }\right)
\end{align*}
This yields
\begin{align*}
   \Esp\left[ \sup_{\xi \in \Xi'_j}  \left(\frac{\left|\hatbeta^a_{j,\xi} -\beta^{a}_{j,\xi}
   \right|}{u_{j,\xi}}\right)^z
\right]
& \leq
\int_{\bR^+} z \tau^{z-1}
\min\left(1,\left|\Xi'_j\right| 2 
e^{-\frac{1}{4} n  \left(\inf_{\xi \in \Xi'_{j}} \frac{u_{j,\xi}}{\sigma_{j,\xi}}  \right)^2\tau^2 }
\right)
d\tau\\
& \quad
+
\int_{\bR^+} z \tau^{z-1}
\min\left(1,\left|\Xi'_j\right| 2 
e^{-\frac{3}{4}n \left(\inf_{\xi \in \Xi'_{j}}
  \frac{u_{j,\xi}}{M_{j,\xi}}  \right) \tau }\right)
d\tau,
\end{align*}
and thus, for any $\tau_1\geq 0$ and $\tau_2 \geq 0$, we get
\begin{align*}
  \Esp\left[ \sup_{\xi \in \Xi'_j}  \left(\frac{\left|\hatbeta^a_{j,\xi} -\beta^{a}_{j,\xi}
   \right|}{u_{j,\xi}}\right)^z
\right]&\leq
\tau_2^z +
\int_{\tau\geq\tau_2} z \tau^{z-1}
\left|\Xi'_j\right| 2
e^{-\frac{1}{4} n \left(\inf_{\xi \in \Xi'_{j}} \frac{u_{j,\xi}}{\sigma_{j,\xi}}\right)^2\tau^2 }
d\tau\\
& \quad+
\tau_1^{z} +
\int_{\tau \geq \tau_1 } z \tau^{z-1}
\left|\Xi'_j\right| 2
e^{-\frac{3}{4}n\left(\inf_{\xi \in \Xi'_{j}} \frac{u_{j,\xi}}{M_{j,\xi}}\right)\tau }
d\tau.
\end{align*}
Take
$$
\tau_1  = \frac{8}{3n} \frac{ \log
  \left(c_{1,z}\left|\Xi'_j\right|\right)}{
\inf_{\xi \in \Xi'_{j}} \frac{u_{j,\xi}}{M_{j,\xi}}}\quad\mathrm{and}\quad
    \tau_2 = \frac{2\sqrt{2}}{\sqrt{n}} \frac{ \sqrt{\log\left(c_{2,z}\left|\Xi'_j\right|\right)}}{
\inf_{\xi \in \Xi'_j} \frac{u_{j,\xi}}{\sigma_{j,\xi}}}.
$$
Hence, by construction, we have:
\begin{align*}
  \forall \tau \geq \tau_1,\ &
\left|\Xi'_j\right| 2
e^{-\frac{3}{4}n \left(\inf_{\xi \in \Xi'_j} \frac{u_{j,\xi}}{M_{j,\xi}}\right) \tau }
\leq \frac{2}{c_{1,z}}
e^{-\frac{3}{8}n\left(\inf_{\xi \in \Xi'_j} \frac{u_{j,\xi}}{M_{j,\xi}}\right)\tau }
\\
  \forall \tau \geq \tau_2,\ &
\left|\Xi'_j\right| 2
e^{-\frac{1}{4} n \left(\inf_{\xi \in \Xi'_{j}}
                               \frac{u_{j,\xi}}{\sigma_{j,\xi}}\right)^2
                               \tau^2 }
\leq \frac{2}{c_{2,z}}
e^{-\frac{1}{8} n \left(\inf_{\xi \in \Xi'_{j}} \frac{u_{j,\xi}}{\sigma_{j,\xi}}\right)^2\tau^2 }.
\end{align*}
This implies
\begin{align*}
& \Esp\left[ \sup_{\xi \in \Xi'_j}  \left(\frac{\left|\hatbeta^a_{j,\xi} -\beta^{a}_{j,\xi}
   \right|}{u_{j,\xi}}\right)^z
\right]\\
& \leq   \left(\frac{2\sqrt{2}}{\sqrt{n}} \frac{ \sqrt{\log\left( c_{2,z}\left|\Xi'_j\right|\right)
}}{\inf_{\xi \in \Xi'_{j}} \frac{u_{j,\xi}}{\sigma_{j,\xi}}} \right)^z
+
2 \left( \frac{2\sqrt{2}}{\sqrt{n}} \frac{1}
{\inf_{\xi \in \Xi'_{j}} \frac{u_{j,\xi}}{\sigma_{j,\xi}}}
\right)^z
\\
&\quad
+
\left(\frac{8}{3n} \frac{ \log \left(c_{1,z}\left|\Xi'_j\right|\right)}
{\inf_{\xi \in \Xi'_j} \frac{u_{j,\xi}}{M_{j,\xi}}}\right)^z
+2
\left(\frac{8}{3n} \frac{1}{\inf_{\xi \in \Xi'_j} \frac{u_{j,\xi}}{M_{j,\xi}}}\right)^z,
\end{align*}
which allows to establish the claimed result.\hfill $\square$\vspace{0.3cm}

Lemma \ref{lem8} allows to obtain the  upper bounds \eqref{elem8i} and \eqref{esuporacle} below.\\
For $u_{j,\xi} = \sigma_{j,\xi}$, we obtain
\begin{align*}
&   \Esp\left[ \sup_{\xi \in \Xi'_j}  \left(\frac{\left|\hatbeta^a_{j,\xi} -\beta^{a}_{j,\xi}
   \right|}{\sigma_{j,\xi}}\right)^z
\right]\\
& \leq   \left(\frac{2\sqrt{2}}{\sqrt{n}} \right)^z  \left(2+
\left(\log \left(c_{2,z}\left|\Xi'_j\right|\right)\right)^{z/2}\right)
+ \left(\frac{8}{3n} \sup_{\xi \in \Xi'_j}
\frac{M_{j,\xi}}{\sigma_{j,\xi}}
\right)^z \left(2+\left(\log\left(c_{1,z}\left|\Xi'_j\right|\right)
\right)^z \right).
\end{align*}
For future use, note that we can also use the uniform bounds $M_j$ (see \eqref{esup}) and
\begin{equation}\label{evarb}
\sigma_{j,\xi}\le \ConstNeedNorm_{2}B(d,2)
B_X^{1/2}2^{j\nu(d)}\myeq \sigma_{j}
\end{equation}
instead of $M_{j,\xi}$ and $\sigma_{j,\xi}$, and obtain
\begin{align}\label{elem8i}
&   \Esp\left[ \sup_{\xi \in \Xi'_j}  \left|\hatbeta^a_{j,\xi} -\beta^{a}_{j,\xi}
   \right|^z
\right]\notag\\
&\leq   \left(\frac{2\sqrt{2}}{\sqrt{n}}\sigma_{j} \right)^z  \left(2+
\left(\log \left(c_{2,z}\left|\Xi'_j\right|  \right)\right)^{z/2} \right)
+ \left(\frac{8}{3n} M_{j} \right)^z \left(2+\left(\log \left(c_{1,z}\left|\Xi'_j\right|
  \right)\right)^z \right).
\end{align}
Along the same lines, with $u_{j,\xi} = c_{\sigma} \sqrt{\log(n)/n}
\sigma_{j,\xi} + c_M \log(n)/ (n-1) M_{j,\xi}$ , we obtain
\begin{align}
 &\Esp \left[
\sup_{\xi \in \Xi'_j}
\left(\frac{
\left|\hatbeta^a_{j,\xi} -\beta^{a}_{j,\xi} \right|
}{
c'_\sigma \sqrt{\log n}
 \frac{\sigma_{j,\xi}}{\sqrt{n}} + c'_M \log n \frac{M_{j,\xi}}{n-1}
}\right)^z
\right]\nonumber
\\
& \leq   \left(\frac{2\sqrt{2}}{c'_\sigma \sqrt{\log n}} \right)^z  \left(2+
\left(\log \left(c_{2,z}\left|\Xi'_j\right|\right)\right)^{z/2}  \right)
+ \left(\frac{8}{3 c'_M \log n}\right)^z \left(2+\left(\log\left(c_{1,z}\left|\Xi'_j\right|\right)\right)^z \right),
\label{esuporacle}
\end{align}
recall that when $\Xi'_j=\Xi_j$, $\left|\Xi'_j\right|\le C_{\Xi}2^{j(d-1)}$.

\subsubsection{Empirical Bernstein and the
probabilities}\label{s722}
We take 
\begin{align*}
\Delta_{j,\xi,\gamma} & =\sqrt{2 \gamma}t_n
\widehat{\sigma}_{j,\xi} + \frac{14}{3}
M_{j,\xi}
\frac{\gamma \log n}{n-1};\\
T_{j,\xi,\gamma} & = 2 \Delta_{j,\xi,\gamma},\quad
T^{b}_{j,\xi,\gamma}= \Delta_{j,\xi,\gamma}, \quad
T^{s}_{j,\xi,\gamma} = 3 \Delta_{j,\xi,\gamma};\\
\Delta^+_{j,\xi,\gamma} &= \sqrt{2 \gamma}t_n
\sigma_{j,\xi} + \frac{26}{3}M_{j,\xi}
\frac{\gamma \log n}{n-1}\quad {\rm and}\quad
\Delta^-_{j,\xi,\gamma} = \sqrt{2 \gamma}t_n
\sigma_{j,\xi} + \frac{2}{3}M_{j,\xi}
\frac{\gamma \log n}{n-1};\\
T^{b,-}_{j,\xi,\gamma} & =  \Delta^{-}_{j,\xi,\gamma}\quad {\rm and}\quad
T^{s,+}_{j,\xi,\gamma} =  3 \Delta^{+}_{j,\xi,\gamma}.
\end{align*}
\begin{lemma}\label{lem9}
The following upper bounds hold:
\begin{align*}
  \Prob&\left\{ T^{b,-}_{j,\xi,\gamma} > T^{b}_{j,\xi,\gamma} \right\}
\leq \frac{1}{n^\gamma};\\
  \Prob&\left\{ T^{s,+}_{j,\xi,\gamma} < T^{s}_{j,\xi,\gamma} \right\}
\leq \frac{1}{n^\gamma};\\
\Prob &\left\{
\left|\hatbeta^a_{j,\xi} -\beta^{a}_{j,\xi}
  \right| > \Delta_{j,\xi,\gamma} \right\} \leq \frac{3}{n^\gamma};\\
&  \Prob\left(
\bigcup_{\xi \in \Quad_j}\left\{T^{s,+}_{j,\xi,\gamma} <
    T^{s}_{j,\xi,\gamma} \right\}\right)
 \leq \sum_{\xi\in\xi_j} \Prob\left\{\ T^{s,+}_{j,\xi,\gamma} <
    T^{s}_{j,\xi,\gamma} \right\} \leq C_{\Xi} 2^{j(d-1)}
  \frac{1}{n^\gamma};\\
&\Prob\left(
\bigcup_{\xi \in \Quad_j}\left\{T^{b,-}_{j,\xi,\gamma}>
  T^{b}_{j,\xi,\gamma}\right\}\right)
 \leq \sum_{\xi\in\xi_j} \Prob\left\{T^{b,-}_{j,\xi,\gamma}>
  T^{b}_{j,\xi,\gamma}\right\}
\leq C_{\Xi} 2^{j(d-1)}
  \frac{1}{n^\gamma};
\\
&\Prob\left(
\bigcup_{\xi \in \Quad_j}\left\{
  \left|\hatbeta^a_{j,\xi}-\beta^{a}_{j,\xi}\right|>\Delta_{j,\xi,\gamma}\right\}\right)
 \leq \sum_{\xi\in\xi_j} \Prob\left\{
  \left|\hatbeta^a_{j,\xi}-\beta^{a}_{j,\xi}\right|>\Delta_{j,\xi,\gamma}\right\}
\leq C_{\Xi} 2^{j(d-1)} \frac{3}{n^\gamma}.
\end{align*}
\end{lemma}
\noindent {\bf Proof.} Using the results of \cite{MP}, we get:
    \begin{align*}
\Prob &\left\{
\sigma_{j,\xi} > \widehat{\sigma}_{j,\xi} + 2\sqrt{2u} \frac{M_{j,\xi}}{\sqrt{n-1}}
\right\} \leq e^{-u};\\
\Prob &\left\{
\sigma_{j,\xi} < \widehat{\sigma}_{j,\xi} - 2\sqrt{2u} \frac{M_{j,\xi}}{\sqrt{n-1}}
\right\} \leq e^{-u};\\
\Prob &\left\{
\left|\hatbeta^a_{j,\xi} -\beta^{a}_{j,\xi}
  \right| > \sqrt{2 u} \frac{\widehat{\sigma}_{j,\xi}}{\sqrt{n}} + \frac{14}{3}
M_{j,\xi}
\frac{u}{n-1}
\right\} \leq 3e^{-u},
\end{align*}
which yields the first inequalities. The others
follow from the union bound.\hfill$\square$\vspace{0.3cm}

\subsubsection{The case $p=\infty$} Let us consider the various terms one by one.\\
{\em Error in the high dimensional
space}.
\begin{align*}
\Esp \left[ M_{j,z} \right] & \leq \Esp \left[ M_{j,z}^{S0} \right]+ \Esp \left[ M_{j,z}^{S1}\right]+ \Esp \left[
      M_{j,z}^{B1}\right]+ \Esp \left[ M_{j,z}^{B2}\right],
\end{align*}
with
\begin{align*}
\Esp \left[ M_{j,z}^{S0} \right ] & =
\sup_{\xi\in\Quad_j}
  \left|
    \beta^{a}_{j,\xi}\right|^z
\ind_{\left|\beta^{a}_{j,\xi}\right|\leq
    T^{s,+}_{j,\xi,\gamma}};\\
\Esp\left[ M_{j,z}^{S1} \right] &
\leq
C_{\Xi} 2^{j(d-1)} \frac{4}{n^\gamma}
\sup_{\xi\in\Quad_j}
  \left|
    \beta^{a}_{j,\xi}\right|^z;
\\
\Esp\left[ M_{j,z}^{B1} \right] & \leq
\Esp\left[ \sup_{\xi \in \Quad_j} \left|\hatbeta^a_{j,\xi} -
    \beta^{a}_{j,\xi}\right|^z
 \ind_{\left|\beta^{a}_{j,\xi}\right|> T^{b,-}_{j,\xi,\gamma}} \right];
\\
\Esp\left[ M_{j,z}^{B2} \right] & \leq
\left(
C_{\Xi} 2^{j(d-1)} \frac{4}{n^\gamma}
\right)^{1-1/\tau}
\left(
\left(\frac{2\sqrt{2}}{\sqrt{n}}\sigma_{j} \right)^{z}  \left(2^{1/\tau}+
\left(\sqrt{\log \left(\left|\Xi_j\right|
      c_{1,{z\tau}}\right)}\right)^{z}  \right)\right.\\
&\mspace{180mu}
\left.
+ \left(\frac{8}{3n} M_{j} \right)^{z} \left(2^{
1/\tau} +\left(\log \left(\left|\Xi_j\right|
  c_{1,{z\tau}}\right)\right)^{z}
\right)\right),
\end{align*}
where we have used $(a+b)^{1/\tau} \leq a^{1/\tau} +
b^{1/\tau}$ for $\tau\geq 1$.\\
This yields
\begin{align*}
&\frac{ \Esp \left[ \left\| \hatfbetam^{a,\rho} -\fbetam^{a,J}
   \right\|_{\infty}^z \right]}{(J+1)^{z-1} \ConstNeedEquiv^z_{\infty}  }\\
& \leq \sum_{j=0}^J
2^{j(d-1)z/2} \left(\sup_{\xi\in\Quad_j}
  \left|
    \beta^{a}_{j,\xi}\right|^z
\ind_{\left|\beta^{a}_{j,\xi}\right|\leq
    T^{s,+}_{j,\xi,\gamma}}+\Esp\left[
\sup_{\xi \in \Quad_j} \left|\hatbeta^a_{j,\xi} -
    \beta^{a}_{j,\xi}\right|^z
 \ind_{\left|\beta^{a}_{j,\xi}\right|> T^{b,-}_{j,\xi,\gamma}}\right]
\right)\\
&
\quad + \frac{4}{n^\gamma}
C_{\Xi}  \sum_{j=0}^J
2^{j(d-1)(z/2+1)}
\sup_{\xi\in\Quad_j}
  \left|
    \beta^{a}_{j,\xi}\right|^z
\\
&\quad
+
\left(
C_{\Xi} \frac{4}{n^\gamma}
\right)^{1-1/\tau}
\sum_{j=0}^J
2^{j(d-1)(z/2+1-1/\tau)}\\
&\quad\quad\quad\quad \times
\left(
\left(\frac{2\sqrt{2}}{\sqrt{n}}\sigma_{j} \right)^{z}  \left(2^{1/\tau}+
\left(\sqrt{\log \left(\left|\Xi_j\right|
      c_{1,{z\tau}}\right)}\right)^{z}  \right)
+ \left(\frac{8}{3n} M_{j} \right)^{z} \left(2^{
1/\tau} +\left(\log \left(\left|\Xi_j\right|
  c_{1,{z\tau}}\right)\right)^{z}
\right)\right)\\
& \myeq O'_{\infty,z} + R'_{1,\infty,z}+R'_{2,\infty,z}.
\end{align*}

\noindent{\em The  terms $R'_{1,\infty,z}$ and $R'_{2,\infty,z}$.} 
The term $R'_{1,\infty,z}$ is the term which appears in Theorem~\ref{toracle}
and thus we only need to bound $R'_{2,\infty,z}$. 
As in the case $p<\infty$, we can use the uniform bounds
  on $\sigma_{j,\xi}$  and $M_{j,\xi}$, namely, \eqref{esup} and \eqref{evarb}, and
  $|\Quad_j| \leq |\Quad_J|$ to obtain
\begin{align*}
&R'_{2,\infty,z}\\
&  \leq
  \left(
 \frac{4C_{\Xi}}{n^\gamma}
\right)^{1-1/\tau}
\sum_{j=0}^J
2^{j(d-1)(z/2+1-1/\tau)}\\
&\quad \times
\left(
\left(\frac{2\sqrt{2}}{\sqrt{n}}\ConstNeedNorm_{2}B(d,2)
2^{j\nu(d)}B_X^{1/2}
\right)^{z}  \left(2^{1/\tau}+
\left(\log \left(c_{2,{z\tau}}\left|\Xi_J\right|
      \right)\right)^{z/2}  \right)\right.\\
&\qquad
\left.+ \left(\frac{8}{3n} \ConstNeedNorm_{\infty}B(d,\infty)2^{j(\nu(d)+(d-1)/2)}B_X
 \right)^{z} \left(2^{
1/\tau} +\left(\log \left(c_{1,{z\tau}}\left|\Xi_J\right|
  \right)\right)^{z}
\right)\right)\\
& \leq
  \left(
 \frac{4C_{\Xi}}{n^\gamma}
\right)^{1-1/\tau}\left[
\left(\frac{2\sqrt{2}}{\sqrt{n}}\ConstNeedNorm_{2}B(d,2)
B_X^{1/2}
\right)^{z}
\left(2^{1/\tau}+
\left(\log \left(\left|\Xi_J\right|
      c_{2,{z\tau}}\right)\right)^{z/2}  \right)
\sum_{j=0}^J
2^{j\left(\nu(d) z+(d-1)(z/2+1-1/\tau)\right)}\right.\\
& \quad\quad\left. +
 \left(\frac{8}{3n} \ConstNeedNorm_{\infty}B(d,\infty)B_X
 \right)^{z}
 \left(2^{
1/\tau} +\left(\log \left(\left|\Xi_J\right|
  c_{1,{z\tau}}\right)\right)^{z}
\right)
\sum_{j=0}^J
2^{j\left(\nu(d) z + (d-1)(z+1-1/\tau)\right)}\right]
\\
& \leq
  \left(
\frac{4C_{\Xi} }{n^\gamma}
\right)^{1-1/\tau}\left[
\left(\frac{2\sqrt{2}}{\sqrt{n}}\ConstNeedNorm_{2}B(d,2)
B_X^{1/2}
\right)^{z}
\left(2^{1/\tau}+
\left(\log \left(c_{2,{z\tau}}\left|\Xi_J\right|
      \right)\right)^{z/2} \right)
\frac{2^{J\left(\nu(d) z+(d-1)(z/2+1-1/\tau)\right)}}{1-2^{-\left(\nu(d) z+(d-1)(z/2+1-1/\tau)\right)}}\right.\\
& \quad\quad\left. +
\left(\frac{8}{3n} \ConstNeedNorm_{\infty}B(d,\infty)B_X
 \right)^{z}
 \left(2^{
1/\tau} +\left(\log \left(c_{1,{z\tau}}\left|\Xi_J\right|
  \right)\right)^{z}
\right)
\frac{2^{J\left(\nu(d) z + (d-1)(z+1-1/\tau)\right)}}
{1-2^{-\left(\nu(d) z + (d-1)(z+1-1/\tau)\right)}}\right].
\end{align*}

\noindent{\em The term $O'_{\infty,z}$.} 
Denote by
\begin{align*}
O_{z,j}'=\sup_{\xi\in\Quad_j}
  \left|
    \beta^{a}_{j,\xi}\right|^z
\ind_{\left|\beta^{a}_{j,\xi}\right|\leq
    T^{s,+}_{j,\xi,\gamma}}+\Esp\left[
\sup_{\xi \in \Quad_j} \left|\hatbeta^a_{j,\xi} -
    \beta^{a}_{j,\xi}\right|^z
 \ind_{\left|\beta^{a}_{j,\xi}\right|> T^{b,-}_{j,\xi,\gamma}}\right].
\end{align*}
Because $T_{j,\xi,\gamma}^{s,++}\ge
T_{j,\xi,\gamma}^{s,+}$, we get
\begin{align*}
&\Esp\left[
\sup_{\xi \in \Quad_j} \left|\hatbeta^a_{j,\xi} -
    \beta^{a}_{j,\xi}\right|^z
 \ind_{\left|\beta^{a}_{j,\xi}\right|> T^{b,-}_{j,\xi,\gamma}}\right]\\
& =  \Esp\left[\sup_{\xi \in \Quad_j} \left|\hatbeta^a_{j,\xi}-\beta^{a}_{j,\xi}\right|^z
\ind_{\left|\beta^{a}_{j,\xi}\right|> T^{s,++}_{j,\xi,\gamma}}\right]+ \Esp\left[\sup_{\xi \in \Quad_j}\left|\hatbeta^a_{j,\xi}-\beta^{a}_{j,\xi}\right|^z
\ind_{T^{s,++}_{j,\xi,\gamma} \geq \left|\beta^{a}_{j,\xi}\right|>
  T^{b,-}_{j,\xi,\gamma}}\right]\\
& \leq  \Esp\left[\sup_{\xi \in \Quad_j} \left|\hatbeta^a_{j,\xi}-\beta^{a}_{j,\xi}\right|^z
\ind_{\left|\beta^{a}_{j,\xi}\right|> T^{s,++}_{j,\xi,\gamma}}\right]\\
&\quad+\Esp\left[ \sup_{\xi \in \Quad_j}\left(\frac{
    \left|\hatbeta^a_{j,\xi}-\beta^{a}_{j,\xi}\right|}{
T^{b,-}_{j,\xi,\gamma}}\ind_{T^{s,++}_{j,\xi,\gamma} \geq \left|\beta^{a}_{j,\xi}\right|>
  T^{b,-}_{j,\xi,\gamma}}\right)^z\right] \sup_{\xi \in \Quad_j}\left\{ \left|\beta^{a}_{j,\xi}\right|^z
\ind_{T^{s,++}_{j,\xi,\gamma} \geq \left|\beta^{a}_{j,\xi}\right|>
  T^{b,-}_{j,\xi,\gamma}}\right\},
\end{align*}
thus
\begin{align*}
  O_{z,j}' & \leq \left(1+\Esp\left[ \sup_{\xi \in \Quad_j}\left(\frac{
    \left|\hatbeta^a_{j,\xi}-\beta^{a}_{j,\xi}\right|}{
T^{b,-}_{j,\xi,\gamma}}\right)^z\right]\right) \sup_{\xi \in \Quad_j}\left\{ \left|\beta^{a}_{j,\xi}\right|^z
\ind_{\left|\beta^{a}_{j,\xi}\right| \leq T^{s,++}_{j,\xi,\gamma} }\right\}\\
&\quad+
\Esp\left[\sup_{\xi \in \Quad_j} \left|\hatbeta^a_{j,\xi}-\beta^{a}_{j,\xi}\right|^z
\ind_{\left|\beta^{a}_{j,\xi}\right|> T^{s,++}_{j,\xi,\gamma}}\right].
\end{align*}
Using now \eqref{esuporacle}, with $c'_{\sigma}=\sqrt{2\gamma}$ and
$c'_M=\frac{2}{3}\gamma$, and $|\Xi_j|\leq C_\Xi 2^{j(d-1)}$,  we get the upper bound in Theorem \ref{toracle}.

\subsubsection{The case $p<\infty$} Let us consider the various terms one by one.\\
\noindent {\em Error in the high dimensional space.} 
We obtain
\begin{align*}
\Esp \left[ S_{j,z} \right] & = \Esp \left[ S_{j,z}^{S0} \right]+ \Esp \left[ S_{j,z}^{S1}\right]+ \Esp \left[
      S_{j,z}^{B1}\right]+ \Esp \left[ S_{j,z}^{B2}\right].
\end{align*}
with
\begin{align*}
  \Esp\left[ S_{j,z}^{S0} \right] & = \sum_{\xi\in\Quad_j}
  \left|
    \beta^{a}_{j,\xi}\right|^z
\ind_{\left|\beta^{a}_{j,\xi}\right|\leq
    T^{s,+}_{j,\xi,\gamma}};\\
  \Esp\left[ S_{j,z}^{S1} \right] &\leq \frac{4}{n^\gamma}  \sum_{\xi\in\Quad_j}
  \left|
    \beta^{
    a}_{j,\xi}\right|^z ;\\
  \Esp\left[ S_{j,z}^{B1} \right] &
 \leq \sum_{\xi\in\Quad_j}
\Esp\left[
  \left|\hatbeta^a_{j,\xi}-\beta^{a}_{j,\xi}\right|^z\right]
\ind_{\left|\beta^{a}_{j,\xi}\right|> T^{b,-}_{j,\xi,\gamma}};\\
  \Esp\left[ S_{j,z}^{B2} \right]
& \leq \frac{4^{1-1/\tau}}{n^{\gamma(1-1/\tau)}} \sum_{\xi\in\Quad_j}
2^{1/\tau} \left( \left( 2 c_{2,z\tau}^{1/(z\tau)}
  \frac{\sigma_{j,\xi}}{\sqrt{n}}\right)^z + \left(\frac{4}{3}
  c_{1,z\tau}^{1/(z\tau)} \frac{M_{j,\xi}}{n} \right)^z\right),
\end{align*}
where we have used $(a+b)^{1/\tau} \leq \left(a^{1/\tau} +
b^{1/\tau}\right)$.  This yields
\begin{align*}
&\frac{\Esp \left[ \left\| \hatfbetam^{a,\rho} -\fbetam^{a,J}
   \right\|_p^z \right]}{(J+1)^{z-1}\ConstNeedEquiv^z_p C_{\Xi}^{z/(p\wedge z)-1}}\\
&\leq \sum_ {j=0}^J
 2^{j(d-1)z(1/2-1/(p\vee z))} \Esp \left[ S_{j,z} \right]\\
& \leq  \sum_ {j=0}^J
 2^{j(d-1)z(1/2-1/(p\vee z))}\sum_{\xi\in\Quad_j} \left(
  \left|
    \beta^{a}_{j,\xi}\right|^z
\ind_{\left|\beta^{a}_{j,\xi}\right|\leq
    T^{s,+}_{j,\xi,\gamma}}
+
\Esp\left[ \left|\hatbeta^a_{j,\xi}-\beta^{a}_{j,\xi}\right|^z\right]
\ind_{\left|\beta^{a}_{j,\xi}\right|> T^{b,-}_{j,\xi,\gamma}}
\right)\\
& \quad + \frac{4}{n^\gamma}  \sum_ {j=0}^J
 2^{j(d-1)z(1/2-1/(p\vee z))} \sum_{\xi\in\Quad_j}
  \left|
    \beta^{
    a}_{j,\xi}\right|^z\\
& \quad + \frac{2^{2-1/\tau}}{n^{\gamma(1-1/\tau)}}\sum_ {j=0}^J
 2^{j(d-1)z(1/2-1/(p\vee z))}
\sum_{\xi\in\Quad_j} \left( \left( 2 c_{2,z\tau}^{1/(z\tau)}
  \frac{\sigma_{j,\xi}}{\sqrt{n}}\right)^z  + \left(\frac{4}{3}
  c_{1,z\tau}^{1/(z\tau)} \frac{M_{j,\xi}}{n} \right)^z \right)\\
& \myeq O_{p,z} + R_{1,p,z} + R_{2,p,z}.
\end{align*}

\noindent{\em The  terms $R_{1,p,z}$ and $R_{2,p,z}$.} 
The term $R_{1,p,z}$ appears as is in Lemma~\ref{toracle}. To bound the term $R_{2,p,z}$, we rely on 
\eqref{esup}. 
We obtain
\begin{align*}
&\sum_{\xi\in\Quad_j} 2^{1/\tau} \left( \left( 2 c_{2,z\tau}^{1/(z\tau)}
  \frac{\sigma_{j}}{\sqrt{n}}\right)^z  + \left(\frac{4}{3}
  c_{1,z\tau}^{1/(z\tau)} \frac{M_{j}}{n} \right)^z \right)
\\
& \leq
\sum_{\xi\in\Quad_j}
2^{1/\tau} \left( 2 c_{2,z\tau}^{1/(z\tau)}
\ConstNeedNorm_{2}B(d,2)
2^{j\nu(d)}B_X^{1/2}
  \frac{1}{\sqrt{n}}\right)^z\\
&\quad
+ \sum_{\xi\in\Quad_j}
2^{1/\tau} \left( \frac{4}{3} c_{1,z\tau}^{1/(z\tau)}
\ConstNeedNorm_{\infty}B(d,\infty)2^{j(\nu(d)+(d-1)/2)}B_X
  \frac{1}{n}\right)^z\\
& \leq
C_{\Xi} 2^{1/\tau} \left( 2 c_{2,z\tau}^{1/(z\tau)}
\ConstNeedNorm_{2}B(d,2) \right)^z B_X^{z/2}
  \frac{1}{n^{z/2}} 2^{j((d-1)+z\nu(d))}
\\ & \quad +
C_{\Xi} 2^{1/\tau} \left( \frac{4}{3} c_{1,z\tau}^{1/(z\tau)}
\ConstNeedNorm_{\infty}B(d,\infty) \right)^z B_X^z
  \frac{1}{n^z} 2^{j((d-1)+z(\nu(d)+(d-1)/2))};
\end{align*}
this yields
\begin{align*}
  &\sum_ {j=0}^J
 2^{j(d-1)z(1/2-1/(p\vee z))} \sum_{\xi\in\Quad_j} 2^{1/\tau} \left( \left( 2 c_{2,z\tau}^{1/(z\tau)}
  \frac{\sigma_{j}}{\sqrt{n}}\right)^z  + \left(\frac{4}{3}
  c_{1,z\tau}^{1/(z\tau)} \frac{M_{j}}{n} \right)^z \right)\\
& \leq \sum_ {j=0}^J
 2^{j(d-1)z(1/2-1/(p\vee z))} C_{\Xi} 2^{1/\tau} \left( 2 c_{2,z\tau}^{1/(z\tau)}
\ConstNeedNorm_{2}B(d,2) \right)^z B_X^{z/2}
  \frac{1}{n^{z/2}} 2^{j((d-1)+z\nu(d))}
\\
& \quad +\sum_ {j=0}^J
 2^{j(d-1)z(1/2-1/(p\vee z))}
C_{\Xi} 2^{1/\tau} \left( \frac{4}{3} c_{1,z\tau}^{1/(z\tau)}
\ConstNeedNorm_{\infty}B(d,\infty) \right)^z B_X^z
  \frac{1}{n^z} 2^{j((d-1)+z(\nu(d)+(d-1)/2))}\\
& \leq C_{\Xi} 2^{1/\tau} \left( 2 c_{2,z\tau}^{1/(z\tau)}
\ConstNeedNorm_{2}B(d,2) \right)^z B_X^{z/2}
  \frac{1}{n^{z/2}} \sum_ {j=0}^J
 2^{jz(\nu(d)+(d-1)/z+(d-1)(1/2-1/(p\vee z)))}
\\ & \quad + C_{\Xi} 2^{1/\tau} \left( \frac{4}{3} c_{1,z\tau}^{1/(z\tau)}
\ConstNeedNorm_{\infty}B(d,\infty) \right)^z B_X^z
  \frac{1}{n^z} \sum_ {j=0}^J
 2^{jz(\nu(d)+(d-1)/z+(d-1)(1-1/(p\vee z))}\\
& \leq \frac{C_{\Xi} 2^{1/\tau} \left( 2 c_{2,z\tau}^{1/(z\tau)}
\ConstNeedNorm_{2}B(d,2) \right)^z}{1-2^{-z(\nu(d)+(d-1)/z+(d-1)(1/2-1/(p\vee z)))}}B_X^{z/2}
  \frac{1}{n^{z/2}} 2^{Jz(\nu(d)+(d-1)/z+(d-1)(1/2-1/(p\vee z)))}
\\ & \quad + \frac{C_{\Xi} 2^{1/\tau} \left( \frac{4}{3} c_{1,z\tau}^{1/(z\tau)}
\ConstNeedNorm_{\infty}B(d,\infty) \right)^z}
{1-2^{-z(\nu(d)+(d-1)/z+(d-1)(1-1/(p\vee z))}}
 B_X^z
  \frac{1}{n^z} 2^{Jz(\nu(d)+(d-1)/z+(d-1)(1-1/(p\vee z))}.
\end{align*}

\noindent {\em The term $O_{p,z}$.} 
Denote by
\begin{align*}
O_{z,j,\xi}=\left|
    \beta^{a}_{j,\xi}\right|^z
\ind_{\left|\beta^{a}_{j,\xi}\right|\leq
    T^{s,+}_{j,\xi,\gamma}}
+
\Esp\left[ \left|\hatbeta^a_{j,\xi}-\beta^{a}_{j,\xi}\right|^z\right]
\ind_{\left|\beta^{a}_{j,\xi}\right|> T^{b,-}_{j,\xi,\gamma}}.
\end{align*}
Because $T_{j,\xi,\gamma}^{s,++}\ge
T_{j,\xi,\gamma}^{s,+}$, we get
\begin{align*}
&\Esp\left[ \left|\hatbeta^a_{j,\xi}-\beta^{a}_{j,\xi}\right|^z\right]
\ind_{\left|\beta^{a}_{j,\xi}\right|> T^{b,-}_{j,\xi,\gamma}}\\
& =  \Esp\left[ \left|\hatbeta^a_{j,\xi}-\beta^{a}_{j,\xi}\right|^z\right]
\ind_{\left|\beta^{a}_{j,\xi}\right|> T^{s,++}_{j,\xi,\gamma}}
+ \Esp\left[ \left|\hatbeta^a_{j,\xi}-\beta^{a}_{j,\xi}\right|^z\right]
\ind_{T^{s,++}_{j,\xi,\gamma} \geq \left|\beta^{a}_{j,\xi}\right|>
  T^{b,-}_{j,\xi,\gamma}}\\
& \leq  \Esp\left[ \left|\hatbeta^a_{j,\xi}-\beta^{a}_{j,\xi}\right|^z\right]
\ind_{\left|\beta^{a}_{j,\xi}\right|> T^{s,++}_{j,\xi,\gamma}}
+ \frac{\Esp\left[
    \left|\hatbeta^a_{j,\xi}-\beta^{a}_{j,\xi}\right|^z\right]}{
\left(T^{b,-}_{j,\xi,\gamma}\right)^z} \left|\beta^{a}_{j,\xi}\right|^z
\ind_{T^{s,++}_{j,\xi,\gamma} \geq \left|\beta^{a}_{j,\xi}\right|>
  T^{b,-}_{j,\xi,\gamma}},
\end{align*}
\begin{align*}
  O_{z,j,\xi} & \leq \left(1+\frac{\Esp\left[
    \left|\hatbeta^a_{j,\xi}-\beta^{a}_{j,\xi}\right|^z\right]}{
\left(T^{b,-}_{j,\xi,\gamma}\right)^z}\right) \left|\beta^{a}_{j,\xi}\right|^z
\ind_{\left|\beta^{a}_{j,\xi}\right| \leq T^{s,++}_{j,\xi,\gamma} } +
\Esp\left[ \left|\hatbeta^a_{j,\xi}-\beta^{a}_{j,\xi}\right|^z\right]
\ind_{\left|\beta^{a}_{j,\xi}\right|> T^{s,++}_{j,\xi,\gamma}}.
\end{align*}
Now using the results of Section~\ref{s721}, with
 $T^{b,-}_{j,\xi,\gamma} = \sqrt{2\gamma}t_n \sigma_{j,\xi} + \frac{2}{3}\gamma \frac{\log n}{n-1}M_{j,\xi}$, we obtain
\begin{align*}
  \sup_{j,\xi} \frac{\Esp\left[
    \left|\hatbeta^a_{j,\xi}-\beta^{a}_{j,\xi}\right|^z\right]}{
\left(T^{b,-}_{j,\xi,\gamma}\right)^z} &\leq 
2 \left(
c_{2,z}
\left(
2 \frac{1}{\sqrt{2\gamma} \sqrt{\log n}}
\right)^z
+
c_{1,z}
\left(
\frac{4}{3}
\frac{1}{(2/3) \gamma \log n  }
\right)^z
\right)\\
& \leq 2 \left(
\left(\frac{\sqrt{2}c_{2,z}^{1/z}}{\sqrt{\gamma \log n}}
\right)^{z}
+
\left(
\frac{2c_{1,z}^{1/z}}{ \gamma \log n  }
\right)^z
\right).
\end{align*}
This yields
\begin{align*}
  O_{p,z} & \leq \left(1+2 \left( \left(
\frac{\sqrt{2} c_{2,z}^{1/z}}{\sqrt{\gamma\log n}}
\right)^z
+
\left(
\frac{2 c_{1,z}^{1/z}}{\gamma \log n}
\right)^z
\right)
\right)\sum_{j=0}^J
 2^{j(d-1)z(1/2-1/(p\vee z))} \\
&\quad\sum_{\xi\in\Quad_j}\left(
\left|\beta^{a}_{j,\xi}\right|^z
\ind_{\left|\beta^{a}_{j,\xi}\right| \leq T^{s,++}_{j,\xi,\gamma} } + \Esp\left[ \left|\hatbeta^a_{j,\xi}-\beta^{a}_{j,\xi}\right|^z\right]
\ind_{\left|\beta^{a}_{j,\xi}\right|> T^{s,++}_{j,\xi,\gamma}}\right).
\end{align*}

\subsection{Proof of Theorem \ref{t2}}\label{s8}
This proof requires an upper bound on: the approximation
error, $R_{1,p,z}$, $R_{1,p,z}$, and $O_{p,z}$. We use that because 
$\fbeta\in B_{r,q}^s(M)$, we have, by Lemma \ref{compnorms}, $\fbetam\in B_{r,q}^s(c_{\mathrm{eq}}M)$.

\subsubsection{The case $1\le p< \infty$}\label{s82} Let us consider the terms one by one.\\
{\em The approximation error.} 
Start with
\[
\left\|\fbetam^{a,J}
  -\fbetam\right\|_p=\left\|\sum_{j>J}\sum_{\xi\in\Xi_j}\beta^a_{j,\xi}\psi_{j,\xi}\right\|_p.
\]
From Lemma \ref{lem1} \eqref{lem1i} and the definition of the Besov
spaces as a sequence space, with $1/q+1/\tilde{q}=1$, we obtain
\[
\left\|\sum_{j>J}\sum_{\xi\in\Xi_j}\beta^a_{j,\xi}\psi_{j,\xi}\right\|_p\le \ConstNeedEquiv_p\sum_{j>J}
2^{-js}2^{j(s+(d-1)(1/2-1/p))}\left\|\left(\beta^a_{j,\xi}\right)_{\xi\in\Xi_j}\right\|_{\ell^p},\]
which yields
\begin{align*}
\left\|\sum_{j>J}\sum_{\xi\in\Xi_j}\beta^a_{j,\xi}\psi_{j,\xi}\right\|_p&\le \ConstNeedEquiv_p2^{-Js}
(2^{s\tilde{q}}-1)^{-1/\tilde{q}}\left\|\fbetam\right\|_{B^s_{p,q}}\\
&\le\left\{\begin{array}{ll}
                                                                       \ConstNeedEquiv_p c_{\mathrm{eq}}M C_{\Xi}^{1/p-1/r}(2^{s\tilde{q}}-1)^{-1/\tilde{q}}\left\|\fbetam\right\|_{B^s_{p,q}}2^{-Js} &\mbox{if $r\ge p$}\\
                                                                       \ConstNeedEquiv_p c_{\mathrm{eq}}M (2^{s\tilde{q}}-1)^{-1/\tilde{q}}\left\|\fbetam\right\|_{B^s_{p,q}}2^{-J(s-(d-1)(1/r-1/p))}&\mbox{if $r\le p$.}
                                                                    \end{array}
\right.
\end{align*}
It is enough to consider the worst case where $r\le p$ and
to check that $\frac{s-(d-1)(1/r-1/p)}{\nu(d)+(d-1)/2}\ge\mu$ in the two zones.\\
In the dense zone, we have 
$$s+\nu(d)+\frac{d-1}{2}\ge\left(\nu(d)+\frac{d-1}{2}\right)\frac{p}{r},$$
which yields $$\frac{s}{s+\nu(d)+\frac{d-1}{2}} \le
\frac{s}{\left(\nu(d)+\frac{d-1}{2}\right)\frac{p}{r}}.$$ Because
$s>(d-1)/r$ and $p\ge r$, we have
$$s-\frac{d-1}{r}+\frac{d-1}{p}-\frac{sr}{p}=(d-1)\left(\frac{sr}{d-1}-1\right)\left(\frac{1}{r}-\frac{1}{p}\right)\ge0,$$
which yields $s-(d-1)(1/r-1/p)\ge \frac{sr}{p}$ and gives the result.\\
In the sparse zone, because $s>(d-1)/r$, we have
$$\frac{s-(d-1)(1/r-1/p)}{\nu(d)+(d-1)/2}\ge\frac{s-(d-1)(1/r-1/p)}{s+\nu(d)-(d-1)(1/r-1/2)}.$$

\noindent{\em The  terms $R_{1,p,p}$ and $R_{2,p,p}$.} 
Using Lemma
\ref{lem2} \eqref{lem2iii} we obtain
\[
R_{1,p,p}\le\frac{4}{n^{\gamma}}(c_{\mathrm{eq}}M)^p\ConstQuadNb^{1-(p\wedge
  r)/r}\sum_{j=0}^J2^{-jp\left(s+(d-1)(1/p-1/(p\wedge r))\right)},
\]
where the exponent is nonpositive because $s>(d-1)/r$, thus
\[
R_{1,p,p}\le\frac{4(c_{\mathrm{eq}}M)^p\ConstQuadNb^{1-(p\wedge
  r)/r}}{n^{\gamma}\left(1-2^{-p\left(s+(d-1)(1/p-1/(p\wedge r))\right)}\right)}.
\]
With $\gamma>p/2$, $R_{1,p,p}$ is of lower order than $t_n^{p}$.\\
We also have
$$R_{2,p,p}\le
\frac{2^{2-1/\tau}}{n^{\gamma(1-1/\tau)}}C_{\Xi}b_{n,p,p,J,\tau}.$$
With the aforementioned choice of $J$, 
\begin{align*}
\frac{1}{\sqrt{n}}2^{J(\nu(d)+(d-1)/2)}B_X^{1/2}&\lesssim
1;\\
\frac{2^{J(d-1)}}{n}B_X&\lesssim
1.
\end{align*}  
Together, these yield that $b_{n,p,p,J,\tau}$ is of the order of a constant.\\
This term is also of lower order than $t_n^{p}$ for $\tau$ large enough such that 
$\gamma(1-1/\tau)>p/2$.\vspace{0.3cm}

\noindent{\em The term $O_{p,p}$.} 
First note that $a_{n,p,p,J}=1+o(1)$.\\
We take $T^{s,++}_{j,\xi,\gamma}$ uniform in $\xi$:
\begin{align*}
  T^{s,++}_{j,\xi,\gamma} &=
 3\sqrt{2 \gamma}t_n
\ConstNeedNorm_{2}B(d,2)2^{j\nu(d)}B_X^{1/2}\\
&\quad+ 52
\ConstNeedNorm_{\infty}B(d,\infty)2^{j(\nu(d)+(d-1)/2)}B_X
\frac{\gamma \log n}{n-1}\\
& \le 2^{j\nu(d)} \sqrt{\gamma} t_n B_X^{1/2}
\left(
3\sqrt{2}\ConstNeedNorm_{2}B(d,2)+
52\ConstNeedNorm_{\infty}B(d,\infty)
  \frac{n\sqrt{\gamma}}{n-1}
\right),
\end{align*}
where the last display uses the upper bound on $J$, this yields,
for $n\ge 2$, 
$$T^{s,++}_{j,\xi,\gamma} 
 \le 2^{j\nu(d)} \sqrt{\gamma} t_n B_X^{1/2}
\left(
3\sqrt{2}\ConstNeedNorm_{2}B(d,2)+
104\ConstNeedNorm_{\infty}B(d,\infty)
\right)\myeq T^{s,++}_{j,\gamma}.
$$
As a consequence of Lemma~\ref{lem7}, we get
\begin{align*}
  \Esp\left[
    \left|\hatbeta^a_{j,\xi}-\beta^{a}_{j,\xi}\right|^p\right]
& \leq 2 \left( \left( 2 c_{2,p}^{1/p}
  \frac{\sigma_{j}}{\sqrt{n}} \right)^p
+ \left( \frac{4}{3} c_{1,p}^{1/p} \frac{M_{j}}{n} \right)^p
\right)\\
& \leq 2 \left( 2 c_{2,p}^{1/p}
\ConstNeedNorm_{2}B(d,2)
2^{j\nu(d)}B_X^{1/2}
  \frac{1}{\sqrt{n}}\right)^p\\
&\quad
+
2\left( \frac{8}{3} c_{1,p}^{1/p}
\ConstNeedNorm_{\infty}B(d,\infty)2^{j(\nu(d)+(d-1)/2)}B_X
  \frac{1}{n}
\right)^p\\
& \leq 2^{jp\nu(d)} \frac{1}{n^{p/2}} B_X^{p/2} 2^{p+1} \left( c_{2,p}^{1/p}
\ConstNeedNorm_{2}B(d,2)+  \frac{4}{3} c_{1,p}^{1/p}
\ConstNeedNorm_{\infty}B(d,\infty)
\right)^p\\
& \leq \frac{\left(T^{s,++}_{j,\gamma}\right)^p}{(\gamma\log n)^{p/2}}
2 \left( 2 \frac{c_{2,p}^{1/p}
\ConstNeedNorm_{2}B(d,2) +  \frac{4}{3} c_{1,p}^{1/p}
\ConstNeedNorm_{\infty}B(d,\infty) }{
3\sqrt{2}\ConstNeedNorm_{2}B(d,2)+
104\ConstNeedNorm_{\infty}B(d,\infty) \left(
  \sqrt{\gamma} \right)
} \right)^p\\
& \leq
\frac{\left(T^{s,++}_{j,\gamma}\right)^p}{(\gamma\log n)^{p/2}}
2 \left( \frac{\sqrt{2}}{3}c_{2,p}^{1/p} +
  \frac{c_{1,p}^{1/p}}{78\sqrt{\gamma}} \right)^p.
\end{align*}

$$\mathrm{Let}\quad 
C_{\gamma} =
3\sqrt{2}\ConstNeedNorm_{2}B(d,2)
+
104
\ConstNeedNorm_{\infty}B(d,\infty) \sqrt{\gamma} \quad\mathrm{and}\quad
\ConstTsigmap=
2^{1/p} \left( \frac{\sqrt{2}}{3}c_{2,p}^{1/p} +
  \frac{c_{1,p}^{1/p}}{78\sqrt{\gamma}} \right).\ \ 
$$
For any $0< z <p$, we have
\begin{align*}
\sum_{\xi\in\Quad_j}& \left(
\left|\beta^{a}_{j,\xi}\right|^p
\ind_{\left|\beta^{a}_{j,\xi}\right| \leq T^{s,++}_{j,\gamma} } + \Esp\left[ \left|\hatbeta^a_{j,\xi}-\beta^{a}_{j,\xi}\right|^p\right]
\ind_{\left|\beta^{a}_{j,\xi}\right|> T^{s,++}_{j,\gamma}}\right)\\
&  \leq \sum_{\xi\in\Quad_j} \left(
\left|\beta^{a}_{j,\xi}\right|^p
\ind_{\left|\beta^{a}_{j,\xi}\right| \leq T^{s,++}_{j,\gamma} } +
\frac{\left(T^{s,++}_{j,\gamma}\right)^p}{(\gamma\log n)^{p/2}}
\ConstTsigmap^p
\ind_{\left|\beta^{a}_{j,\xi}\right|> T^{s,++}_{j,\xi,\gamma}}
\right)\\
& \leq \left(1+\frac{\ConstTsigmap^p}{(\gamma\log n)^{p/2}}\right)\left( T^{s,++}_{j,\gamma}\right)^{p-z}
\sum_{\xi\in\Quad_j} \left|\beta^{a}_{j,\xi}\right|^{z}
\\
&\leq \left(1+\frac{\ConstTsigmap^p}{(\gamma\log n)^{p/2}}\right)\left(\sqrt{\gamma}t_n
B_X^{1/2}C_{\gamma}
\right)^{p-z} 2^{j\nu(d)(p-z)}\sum_{\xi\in\Quad_j}
\left|\beta^{a}_{j,\xi}\right|^{z}.
\end{align*}
We need to sum over $j$ and take two different values for
$z$, one that we denote $z_1$ for $j\le j_0$ and one that we denote
$z_2$ for $j_0<j\le J$. The values of $z_1$, $z_2$, $j_0$ will be specified later,
depending on the value of the parameters $r,\ q,\ s$ and $p$ such
that we are in the dense or sparse zone.  Up to a multiplying
constant, we thus need to control
\begin{align*}
A+B=&\left(B_X^{1/2}t_n\right)^{p-z_1}\sum_{j=0}^{j_0}2^{j\left[\nu(d)(p-z_1)+(d-1)(p/2-1)\right]}\sum_{\xi\in\Xi_j}\left|\beta_{j,\xi}^a\right|^{z_1}\\
&+\left(B_X^{1/2}t_n\right)^{p-z_2}\sum_{j=j_0+1}^{J}2^{j\left[\nu(d)(p-z_2)+(d-1)(p/2-1)\right]}\sum_{\xi\in\Xi_j}\left|\beta_{j,\xi}^a\right|^{z_2},
\end{align*}
where we choose adequately $z_1$, $z_2$ and $j_0$ in the two zones.  Because of Lemma \ref{lem2} \eqref{lem2i}, we only consider $p\ge r$.\\
Let us first consider the dense zone.  We define
\[
\tilde{r}=\frac{p(\nu(d)+(d-1)/2)}{s+\nu(d)+(d-1)/2}.
\]
In the dense zone, $\tilde{r}\le r$, $p>\tilde{r}$ and
\begin{equation}\label{evaraux}
s=\left(\nu(d)+\frac{d-1}{2}\right)\left(\frac{p}{\tilde{r}}-1\right).
\end{equation}
With $z_2=r$, we get
\[
B\le
\left(B_X^{1/2}t_n\right)^{p-r}\sum_{j=j_0+1}^{J}2^{j\left[\nu(d)(p-r)+(d-1)(p/2-1)\right]}\sum_{\xi\in\Xi_j}\left|\beta_{j,\xi}^a\right|^{r}.
\]
\noindent Lemma \ref{lem2} \eqref{lem2iii} gives that
\[
\sum_{\xi\in\Quad_j}
|\beta_{j,\xi}|^r
\leq D_{j}^{r}
2^{-jr(s+(d-1)(1/2-1/r))},
\]
where $\forall j\in \xN,\ D_j\ge0$, $(D_{j})_{j\in\xN} \in \ell_q$.
Note that
\begin{equation}\label{edec}
s+(d-1)\left(\frac12-\frac{1}{r}\right)=\frac{(d-1)p}{2\tilde{r}}-\frac{d-1}{r}+\nu(d)\left(\frac{p}{\tilde{r}}-1\right),
\end{equation}
thus
\begin{align*}
B&\le \left(B_X^{1/2}t_n\right)^{p-r}\sum_{j=j_0+1}^{J}2^{jp\left(1-\frac{r}{\tilde{r}}\right)\left(\nu(d)+\frac{d-1}{2}\right)}D_{j}^{r}\\
&\lesssim (c_{\mathrm{eq}}M)^r \left(B_X^{1/2}t_n\right)^{p-r}2^{j_0p\left(1-\frac{r}{\tilde{r}}\right)\left(\nu(d)+\frac{d-1}{2}\right)},
\end{align*}
for $q\ge1$ if $r>\tilde{r}$ and for $q\le r$ if $r=\tilde{r}$ (i.e.,
$s=p\left(\nu(d)+\frac{d-1}{2}\right)\left(\frac{1}{r}-\frac{1}{p}\right)$).\\
Taking
$2^{j_0\frac{p}{\tilde{r}}\left(\nu(d)+\frac{d-1}{2}\right)}\simeq
\left(B_X^{1/2}t_n\right)^{-1}$, 
we get
\[
B\lesssim M^r
\left(B_X^{1/2}t_n\right)^{p-\tilde{r}},
\]
which is the rate that we expect in that zone.\\
As for $A$, we take $z_1=\overline{r}<\tilde{r}\le r$, this yields,
using Lemma \ref{lem2} \eqref{lem2iii},
\begin{align*}
A&\le \left(B_X^{1/2}t_n\right)^{p-\overline{r}}\sum_{j=0}^{j_0}2^{j\left[\nu(d)(p-\overline{r})+(d-1)(p/2-1)\right]}\sum_{\xi\in\Xi_j}\left|\beta_{j,\xi}^a\right|^{\overline{r}}\\
&\lesssim M^r  \left(B_X^{1/2}t_n\right)^{p-\overline{r}}\sum_{j=0}^{j_0}2^{j\left[\nu(d)(p-\overline{r})+(d-1)(p/2-1)-\overline{r}\left(s+(d-1)(1/2-1/\overline{r})\right)\right]}\\
&\lesssim M^r  \left(B_X^{1/2}t_n\right)^{p-\overline{r}}\sum_{j=0}^{j_0}2^{jp(\nu(d)+(d-1)/2)(1-\overline{r}/\tilde{r})}\quad \mbox{(using \eqref{evaraux})}\\
&\lesssim M^r  \left(B_X^{1/2}t_n\right)^{p-\overline{r}}2^{j_0p(\nu(d)+(d-1)/2)(1-\overline{r}/\tilde{r})}\\
&\lesssim M^r  \left(B_X^{1/2}t_n\right)^{p-\tilde{r}}\quad \mbox{(from the definition of $j_0$)}.
\end{align*}
Let us now consider the sparse zone.  We define by
\[
\tilde{r}=p\frac{\nu(d)+(d-1)(1/2-1/p)}{s+\nu(d)-(d-1)(1/r-1/2)},
\]
in a such a way that
\begin{align}
p-\tilde{r}&=p\frac{s-(d-1)(1/r-1/p)}{s+\nu(d)-(d-1)(1/r-1/2)};\notag\\
\tilde{r}-r&=\frac{(p-r)((d-1)/2+\nu(d))-rs}{s+\nu(d)-(d-1)(1/r-1/2)}>0;\notag\\
s+(d-1)\left(\frac12-\frac{1}{r}\right)&=\frac{(d-1)p}{2\tilde{r}}-\frac{d-1}{\tilde{r}}+\nu(d)\left(\frac{p}{\tilde{r}}-1\right).\label{edec2}
\end{align}
For the term $A$, we take $z_1=r$ and obtain
\begin{align*}
A&\le\left(B_X^{1/2}t_n\right)^{p-r}\sum_{j=0}^{j_0}2^{j\left[\nu(d)(p-r)+(d-1)(p/2-1)\right]}
\sum_{\xi\in\Xi_j}\left|\beta_{j,\xi}^a\right|^{r}\\
&\le\left(B_X^{1/2}t_n\right)^{p-r}\sum_{j=0}^{j_0}2^{j\left[\nu(d)+(d-1)(1/2-1/p)\frac{p}{\tilde{r}}(\tilde{r}-r)\right]}D_{j}^{r}
\quad \mbox{(using \eqref{edec2})}\\
&\lesssim\left(B_X^{1/2}t_n\right)^{p-r}2^{j_0\left[(\nu(d)+(d-1)(1/2-1/p)\frac{p}{\tilde{r}}(\tilde{r}-r)\right]}M^r,
\end{align*}
the last inequality holds because $\nu(d)+(d-1)/2-(d-1)/p>0$, indeed,
because we are in the sparse zone $\nu(d)+(d-1)/2\ge
s/(p/r-1)=sr/(p-r)\ge 2/(p-r)\ge (d-1)/p$.  Taking
$2^{j_0(\nu(d)+(d-1)(1/2-1/p))\frac{p}{\tilde{r}}}\simeq\left(B_X^{1/2}t_n\right)^{-1}$,
yields 
\[A\lesssim
M^r\left(B_X^{1/2}t_n\right)^{p-\tilde{r}}.
\]
For the term $B$, we take $z_2=\overline{r}>\tilde{r}>r$ and obtain
\begin{align*}
B&\le\left(B_X^{1/2}t_n\right)^{p-\overline{r}}\sum_{j=j_0+1}^{J}2^{j\left[\nu(d)(p-\overline{r})+(d-1)(p/2-1)\right]}
\sum_{\xi\in\Xi_j}\left|\beta_{j,\xi}^a\right|^{\overline{r}}\\
&\lesssim\left(B_X^{1/2}t_n\right)^{p-\overline{r}}\sum_{j=j_0+1}^{J}2^{j(\nu(d)+(d-1)(1/2-1/p))p(r-\overline{r})/\tilde{r}}
D_j^{\overline{r}}\quad \mbox{(using \eqref{edec2})}\\
&\lesssim\left(B_X^{1/2}t_n\right)^{p-\overline{r}}2^{j_0(\nu(d)+(d-1)(1/2-1/p))p(r-\overline{r})/\tilde{r}}
M^{\overline{r}}\\
&\lesssim\left(B_X^{1/2}t_n\right)^{p-\tilde{r}}M^{\overline{r}}.
\end{align*}

\subsubsection{The case $p=\infty$}\label{s81}
Consider $r=\infty$. The general case follows by Lemma \ref{lem2} \eqref{lem2ii}.\\ 
\noindent{\em The approximation error.} 
Because $f_{\beta}\in B^s_{\infty,q}(M)$, we have by 
Lemma \ref{lem1} \eqref{lem1i} 
\begin{align*}
\left\|\sum_{j>J}\sum_{\xi\in\Xi_j}\beta^a_{j,\xi}\psi_{j,\xi}\right\|_{\infty}&\le\sum_{j>J}\left\|\sum_{\xi\in\Xi_j}\beta^a_{j,\xi}\psi_{j,\xi}\right\|_\infty\\
&\le C_{\infty}'c_{\mathrm{eq}}M\sum_{j>J}2^{j(d-1)/2}2^{-j(s+(d-1)/2)}D_j\quad (\mathrm{where\ }\|(D_j)_{j\in\xN}\|_q\le c_{\mathrm{eq}}M)\\
&\le C_{\infty}'c_{\mathrm{eq}}M2^{-Js}(2^{s\tilde{q}}-1)^{-1/\tilde{q}}.
\end{align*}
From the choice of $J$, we get
\[
\left\|\sum_{j>J}\sum_{\xi\in\Xi_j}\beta^a_{j,\xi}\psi_{j,\xi}\right\|_{\infty}\lesssim  C_{\infty}'c_{\mathrm{eq}}M(2^{s\tilde{q}}-1)^{-1/\tilde{q}}
\left(t_nB_X^{1/2}\right)^{s/(\nu(d)+(d-1)/2)}.
\]
This term is negligible because $s/(\nu(d)+(d-1)/2)\ge
s/(s\nu(d)+(d-1)/2)$.\vspace{0.3cm}

\noindent {\em The terms $R'_{1,\infty,z}$ and $R'_{2,\infty,z}$.} 
Using the definition of the Besov norm, 
we obtain
\begin{align*}
R'_{1,\infty,z}&\le\frac{4}{n^{\gamma}}(c_{\mathrm{eq}}M)^z\ConstQuadNb\sum_{j=0}^J2^{-jzs}2^{j(d-1)}\\
&\lesssim\frac{4}{n^{\gamma}}2^{J(d-1)}M^z.
\end{align*}
With $\gamma>z/2+1$, which holds if $2(\gamma-1)(1-1/\tau)>z$, $R_{1,\infty,z}$ is of lower order than $t_n^{z}$.\\
Due to the choice of $J$, the term in bracket in the expression of $R'_{2,\infty,z}$
in Theorem \ref{toracle} is less than 1. The second term in the expression of $b_{n,\infty,z,J,\tau}$ is of smaller order than the first term. The order of $b_{n,\infty,z,J,\tau}$ is finally $(\log n)^{z/2}$.
Thus, we have
$$R'_{2,\infty,z}
\lesssim
\left(n^{-\gamma}2^{J(d-1)}\right)^{1-1/\tau}(\log n)^{z/2}.$$
This term is also of lower order than $t_n^{z}$ when $\tau$ is such that $2(\gamma-1)(1-1/\tau)>z$.\vspace{0.3cm}

\noindent {\em The term $O'_{\infty,z}$.} 
Note that here $a_{n,\infty,z,J}$ is of the order of a constant.  
We now proceed like for the term $O_{p,p}$. 
Using \eqref{elem8i}, we obtain for arbitrary $\overline{z}\in[0,z]$
\begin{align*}
&\sup_{\xi \in \Quad_j}
\left|\beta^{a}_{j,\xi}\right|^z
\ind_{\left|\beta^{a}_{j,\xi}\right| \leq T^{s,++}_{j,\xi,\gamma} } + \Esp\left[ \sup_{\xi \in \Quad_j}\left|\hatbeta^a_{j,\xi}-\beta^{a}_{j,\xi}\right|^z\ind_{\left|\beta^{a}_{j,\xi}\right|> T^{s,++}_{j,\xi,\gamma}}\right]\\
&\lesssim \left(\sqrt{\gamma}t_n
B_X^{1/2}
\right)^{z-\overline{z}} 2^{j\nu(d)(z-\overline{z})}\sup_{\xi\in\Quad_j}
\left|\beta^{a}_{j,\xi}\right|^{\overline{z}}.
\end{align*}
 We use an upper bound on $A+B$, where:
\begin{align*}
A&=\left(B_X^{1/2}t_n\right)^{z-z_1}\sum_{j=0}^{j_0}2^{j\left[\nu(d)(z-z_1)+(d-1)z/2\right]}
\quad\sup_{\xi\in\Xi_j}\left|\beta_{j,\xi}^a\right|^{z_1};\\
B&=\left(B_X^{1/2}t_n\right)^{z-z_2}\sum_{j=j_0+1}^{J}2^{j\left[\nu(d)(z-z_2)+(d-1)z/2\right]}
\sup_{\xi\in\Xi_j}\left|\beta_{j,\xi}^a\right|^{z_2},
\end{align*}
for well-chosen $0\le j_0\le J$, $z_1$ and $z_2$.  Because $f\in B^s_{\infty,q}(M)$, we have
\[\forall  \overline{z}\ge1,\ 
\sup_{\xi\in\Xi_j}\left|\beta_{j,\xi}^a\right|^{\overline{z}}\le
(c_{\mathrm{eq}}M)^{\overline{z}}2^{-j(s+(d-1)/2)\overline{z}}.
\]
The result follows taking $z_1=0$, $j_0$ such that $2^{j_0}\simeq t_n^{-1/(s+\nu(d)+(d-1)/2)}$, and $z_2=z$.

\end{document}